\documentclass[12pt,draft]{amsart}
\usepackage{amsmath,amsthm,amsfonts,amssymb}

\date{}
\newtheorem{theorem}{Theorem }[section]
\newtheorem{lemma}[theorem]{Lemma}
\newtheorem{remark}[theorem]{Remark}%[section]
\newtheorem{corollary}[theorem]{Corollary}%[section]
\newtheorem{example}[theorem]{Example}
\newtheorem{definition}[theorem]{Definition}

\newcommand{\MM}{{\mathcal M}}

\def\R{\mathbb{R}}
\def\R{\mathbb{R}}

\begin{document}

\title{Modified log-Sobolev inequalities and isoperimetry}

\author{Alexander V. Kolesnikov}

\address{
Department of Mechanics and Mathematics,
         Moscow State University, 119992 Moscow, Russia. e-mail: sascha77@mail.ru}

\begin{abstract}
We find sufficient conditions for a probability measure $\mu$ to satisfy
an inequality of the type
$$
\int_{\R^d} f^2 F\Bigl( \frac{f^2}{\int_{\R^d} f^2 \,d \mu} \Bigr) \,d \mu
\le C \int_{\R^d} f^2 c^{*}\Bigl( \frac{|\nabla f|}{|f|} \Bigr) \,d \mu
+ B \int_{\R^d} f^2 \,d \mu,
$$
where $F$ is concave and $c$ (a~cost function) is convex.
We show that under broad assumptions on $c$ and $F$
 the above inequality
 holds if for some $\delta>0$ and $\varepsilon>0$
 one has
 $$
 \int_{0}^{\varepsilon} \Phi\Bigl(\delta c\Bigl[\frac{t F(\frac{1}{t})}{{\mathcal I}_{\mu}(t)} \Bigr] \Bigr) \,dt
 < \infty,
 $$
where ${\mathcal I}_{\mu}$ is the isoperimetric function of $\mu$ and $\Phi = (y F(y) -y)^{*}$.
In a partial case
$${\mathcal I}_{\mu}(t) \ge k t \varphi ^{1-\frac{1}{\alpha}} (1/t),
$$
where $\varphi$ is a concave function growing not faster than $\log$,  $k>0$, $1 < \alpha \le 2$
 and  $t \le 1/2$, we establish a family of
tight inequalities interpolating between the
$F$-Sobolev and modified  inequalities of log-Sobolev type.
A basic example is given by
 convex measures  satisfying certain integrability assumptions.
\end{abstract}

\maketitle

\section{Introduction}

The celebrated logarithmic Sobolev inequality
\begin{equation}
\label{LSI}
\mbox{Ent}_{\mu} f^2:=
\int_{\R^d} f^2 \log \Bigl( \frac{f^2}{\int_{\R^d} f^2 \,d \mu} \Bigr) \,d \mu
\le
2C \int_{\R^d} \bigl| \nabla f \bigr|^2 \,d \mu,
\end{equation}
where $\mu=e^{-V}\,dx$ is a probability measure,
has numerous applications in probability theory, mathematical physics, and geometry.
It appeared first in the work of Gross \cite{Gross}, where he established (\ref{LSI})
 for the standard Gaussian measure.
Gross discovered that (\ref{LSI}) implies  hypercontractivity of the semigroup
$e^{tL}$ generated by
 $ L=\Delta - \bigl<\nabla V, \nabla  \bigr>.$

Necessary and sufficient conditions for (\ref{LSI}) have been intensively
studied by many authors (see \cite{ISLFr}).
It is well-known  that for every probability measure satisfying (\ref{LSI}) there exists
$\varepsilon >0$
such that
\begin{equation}
\label{expbound}
e^{\varepsilon |x|^2} \in L^1(\mu).
\end{equation}
It has been shown by Wang (\cite{Wang97})
that this assumption is sufficient provided $\mu$
is convex, i.e., has the form
$\mu = e^{-V} \,d x$,
 where $V$ is a convex function (in the literature  convex measures are also called log-concave).
 Wang's proof employs the associated diffusion semigroup.
 Bobkov  \cite{Bobkov} gave another proof of this result by applying the Pr{\'e}kopa--Leindler
 theorem and isoperimetric inequalities.
 There exist non-convex measures satisfying   (\ref{LSI}).
For example, according to a result of Holley and Strook, if $\mu$ satisfies
(\ref{LSI}), every probability measure  $e^{\varphi} \cdot \mu$
with $ a \le \varphi  \le b$ satisfies logarithmic Sobolev
inequality with $C' = e^{2(b-a)} C$.

 Recall that (\ref{LSI}) implies the Poincar{\'e} inequality
\begin{equation}
\label{Poin}
\mbox{Var}_{\mu} f:= \int_{\R^d} f^2 \,d \mu - \Bigl(\int_{\R^d} f \,d \mu \Bigr)^2
\le C \int_{\R^d} | \nabla f |^2 \,d \mu.
\end{equation}

The log-Sobolev inequality can be considered as a Poincar{\'e}-type inequality
for the $L^2 \log L$-Orlicz norm. By using this observation  and some classical results
on Hardy's inequality  with weights, Bobkov and G{\"o}tze \cite{BG} established necessary and sufficient conditions for
(\ref{LSI}) on the real line. Namely, $\mu = \rho \,dx$ satisfies (\ref{LSI}) if and only if
$$
\sup_{x<m} F(x) \log \Bigl( \frac{1}{F(x)} \Bigr) \int_{x}^{m} \frac{\,dx}{\rho(x)}
< \infty,
$$
$$
\sup_{x>m} (1-F(x)) \log \Bigl( \frac{1}{1-F(x)} \Bigr) \int_{m}^{x} \frac{\,dx}{\rho(x)} < \infty,
$$
where $F(x) = \mu \Bigl( (-\infty,x] \Bigr)$ and $m$ is the median of $\mu$.

It is well-known that (\ref{LSI}) (as well as the classical Sobolev inequalities)
 is
closely related to the isoperimetric inequalities.
For every Borel $A \subset \R^d$ we denote by $\mu^{+}(A)$ the surface measure of the boundary $\partial A$:
$$
\mu^{+}(A) = \underline{\lim}_{h \to 0} \frac{\mu(A^{h}) - \mu(A)}{h},
$$
where $A^{h} = \{x : \mbox{dist}(x,A) \le h \}$ is the
$h$-neighborhood of $A$. It was proved by Ledoux \cite{Led94}
that the isoperimetric inequality of the Gaussian type
$$
\mu^{+}(A) \ge c \varphi \bigl( \Phi^{-1}(\mu(A)) \bigr)
$$
implies (\ref{LSI}). Here
$$
\varphi(x) = \frac{1}{\sqrt{2 \pi}} e^{-\frac{x^2}{2}}, \ \Phi(x) = \int_{-\infty}^{x} \varphi(s) \,ds.
$$

Some sufficient conditions for (\ref{LSI}) can be obtained by perturbation
methods. For example, Carlen and Loss applied in \cite{CL}
the log-Sobolev inequality
$$
\int_{\R^d} f^2 \log f^2 \,dx \le \frac{1}{\pi e^2}
\int_{\R^d} | \nabla f|^2 \,dx, \ \int_{\R^d} f^2 \,dx =1
$$
for Lebesgue measure. In particular, they
proved  that $\mu = e^{-V} \,dx$ satisfies (\ref{LSI}) provided
that
$$
\frac{1}{4} | \nabla V|^2 - \frac{1}{2} \Delta V - \pi e^2 V
$$
is bounded from below and $\mu$ satisfies  (\ref{Poin})
(see also \cite{BCR} and \cite{Cat}).

It follows from  (\ref{expbound})  that $\mu$ has a very fast decay.
However, many distributions exhibit some weaker, yet useful properties.
Below we consider the following generalizations of (\ref{LSI}):

 1) The defective log-Sobolev inequality
 $$
 \mbox{Ent}_{\mu} f^2 \le
2C \int_{\R^d} \bigl| \nabla f \bigr|^2 \,d \mu
+ B \int_{\R^d} f^2 \,d \mu.
 $$

2) The $F$-Sobolev inequality
$$
\int_{\R^d} f^2 F\Bigl(\frac{f^2}{\int_{\R^d} f^2 \,d\mu }\Bigr) \,d \mu \le
 2C \int_{\R^d} \bigl| \nabla f \bigr|^2 \,d \mu + B \int_{\R^d} f^2 \,d\mu ,
$$
where $F$ is a concave function.

3) The modified log-Sobolev inequality
\begin{equation}
\label{MLSI1}
\mbox{\rm{Ent}}_{\mu} f^2 \le
C \int_{\R^d} f^2 c^{*} \Bigl( \frac{|\nabla f|}{|f|}\Bigr) \,d \mu
\end{equation}
for some convex
$c : \R^{+} \to \R^{+}$. Here
$c^{*}(x) = \sup_{y \in \R^{+}} \bigl( \bigl<x,y\bigr> - c(y)\bigr)$.

Inequality of  type 1) implies the hyperboundedness of the associated semigroups
(see \cite{Dav}).
A basic example for 2) and 3) is given by the following measure on the real line:
$$
\mu_{\alpha} = Z_{\alpha} e^{-|x|^{\alpha}} \,d x,
$$
where $1< \alpha \le 2$. It was proved in \cite{GGM} that $\mu_{\alpha}$
satisfies (\ref{MLSI1}) with
\begin{equation}
\label{cost}
c(x)=c_{A, \alpha}(x) =  \left\{
\begin{array}{lcr}
\frac{x^2}{2}\  \quad       \mbox{if} \ |x| \le A \\
A^{2-\alpha} \frac{|x|^{\alpha}}{\alpha} + A^2 \frac{\alpha-2}{2\alpha} \ \quad \mbox{if} \ |x| \ge A, \\
\end{array}
\right.
\end{equation}
for every $A>0$.  By the tensorization argument
the result holds also in the multidimensional case for the product measure
$\prod_{i=1}^{d} \mu_{\alpha}(\,d x_i)$ and the
cost function $c_{d, A, \alpha}(x) = \sum_{i=1}^{d} c_{A,\alpha}(x_i)$.
 On the other hand, by a result from \cite{BCR},  $\mu_{\alpha}$ satisfies
\begin{align*}
\int f^2 \log^{2/\beta} (1+f^2) \,d \mu - \Bigl(\int f^2 \,d \mu \Bigr)
\log^{2/\beta}\Bigl(1+ \int f^2 \,d \mu\Bigr)
\le C \int |\nabla f|^2 \,d \mu,
\end{align*}
where $\frac{1}{\alpha} + \frac{1}{\beta} =1$.
One can easily verify that $c^{*}_{A,\alpha} = c_{A,\beta}$.

The case $\alpha  \ge 2$ has been considered in  \cite{BZ}.
 In this case the measure
 $$\mu = Z_{\alpha,d} \prod_{i=1}^{d} e^{-|x_i|^{\alpha}}\,dx_i$$
on $\R^d$ satisfies the inequality
\begin{equation}
\label{Bob-Zeg}
\mbox{Ent}_{\mu} |f|^{\beta} \le C \int_{\R^d}
\sum_{i} |\partial_{x_i} f|^{\beta} \,d \mu.
\end{equation}

Among other  generalizations of (\ref{LSI}) let us  mention
an important result from \cite{LO} on a family of  inequalities interpolating
between log-Sobolev and
 Poincar{\'e}. If  $1 <\alpha \le 2$, $1 \le p \le 2$,
 then for every smooth $f$ one has
$$
\int_{\R^d} f^2 \,d \mu_{\alpha} -
\Bigl( \int_{\R^d} |f|^p \,d \mu_{\alpha}\Bigr)^{2/p}
\le
C(2-p)^{2\bigl(1-\frac{1}{\alpha}\bigr)} \int_{\R^d} |\nabla f|^2 \,d \mu_{\alpha}.
$$
Inequalities of this type were proved first by Beckner in \cite{Beck} for Gaussian measures.
For further development  and connections with the $F$-Sobolev inequality,
see  \cite{BCR}, \cite{BCR2}, and \cite{Wang2005}.

Inequality (\ref{MLSI1}) is closely related to  the  Talagrand transportation inequality
\begin{equation}
\label{c-Tal}
W_c(\mu, f \cdot \mu) \le \mbox{Ent}_{\mu} f,
\end{equation}
where $f \cdot \mu$ is another probability measure and $W_c$ is the minimum of the
Kantorovich
functional for  the cost  function $c$ (see \cite{Vill} for details).
In fact, under broad assumptions
on~$c$, inequality (\ref{MLSI1})
is  stronger than (\ref{c-Tal}). This was shown
in \cite{OV} for the case of quadratic cost function.
 It was proved  in \cite{CGH} by  the optimal transportation method  that (\ref{MLSI1})
 holds
 for measures of the type
$\mu = e^{-V}\,dx$, where $V$ satisfies
$$
V(b)-V(a) \ge  \bigl<\nabla V(a), b-a\bigr> + \alpha c(b-a)
$$
for some $\alpha>0$ and a proper choice of $c$. For
recent
progress in transportation inequalities, including some
exponential- and power-type  estimates, see \cite{BK},
\cite{BV}, \cite{Goz1}, \cite{Goz2},  \cite{Kol},
and the references therein.

In this paper we obtain
sufficient conditions which guarantee inequalities
 of the following type:
\begin{equation}
\label{main}
\int_{\R^d} f^2 F\Bigl(\frac{f^2}{\int_{\R^d} f^2 \,d \mu}\Bigr) \,d \mu \le
C \int_{\R^d} f^2 c^{*} \Bigl( \frac{|\nabla f|}{|f|}\Bigr) \,d \mu
+ B  \int_{\R^d} f^2 \,d \mu,
\end{equation}
where $F$ is concave and $c: \R^{+} \to \R^{+}$ is convex (Theorem \ref{mainth}).
This inequality unifies the defective modified log-Sobolev inequalities and the
$F$-Sobolev inequalities.
Obviously, the tight $F$-Sobolev inequality corresponds
to the case $c=|x|^2$, $B=0$,
and the modified Sobolev inequality corresponds to
the case $F=\log$, $B=0$.

 An important assumption on $c$ which we use below (though not everywhere)
 is the following:
 \begin{equation*}\label{homog}
\mbox{for any $k>0$ there is $n(k)>0$ such that } c(kx) \le n(k) c(x), \quad
 c^{*}(kx)  \le n(k) c^{*}(x).
 \leqno (H)
 \end{equation*}

 Our estimate is based on the use of a special isoperimetric function
$$
I_{F}(r) = \sup_{A \in \MM_r} \frac{\mu(A) F \bigl( \frac{1}{\mu(A)} \bigr) }{\mu^{+}(A)}.
$$
Here $\MM_r = \{A : \mu(A) = \mu(\{x: |x| > r\})\}$. Assume that
(H) holds. The main result (Theorem \ref{mainth}, Remark
\ref{scal}) can be roughly formulated in the following way:

\vskip .1in

{\it Integrability of
$\Phi(\delta c(I_{F}))$
for some $\delta>0$, where
$\Phi = (yF(y)-y)^{*}$, implies {\rm(\ref{main})}.
}

\vskip .1in

Let us give some important examples of the function $I_F$.
In the case of a convex measure $\mu$ and  $F =\log $,
the function $I_F(r)$ can be estimated for large values $r$ of  by
$C r$ with some $C>0$. This follows from an estimate obtained  in   \cite{BH}
(see Lemma \ref{conv-iso}).
In the case of an entropy functional $F$ growing
as $\log^{\tau}(x)$, $\tau \le 1$ and additional assumption that $\exp(|x|^{\alpha}) \in L^1(\mu)$, this result
combined with Chebyshev's inequality yields
that $I_{F}(r) \le C r^{1-\alpha(1-\tau)}$
(see Lemma \ref{isoest} for a precise result).

The integrability assumption can be rewritten even in a more elegant way if
we employ the classical isoperimetric function
${\mathcal I}_{\mu}$  of $\mu$ defined by
\begin{equation}
\label{IF}
{\mathcal I}_{\mu}(t) = \inf_{A \subset \R^d: \mu(A)=t} \mu^{+}(A).
\end{equation}
Assume that $c$  satisfies (H).
It turns out that (\ref{main}) holds for a broad class
of $F$ and $c$ if for some $\delta>0$, $K>1$ one has
 \begin{equation}
 \label{12.11}
 \int_{0}^{1/K}
 \Phi\Bigl(\delta c
 \Bigl[\frac{t F(\frac{1}{t})}{{\mathcal I}_{\mu}(t)}
 \Bigr] \Bigr) \,dt;
 < \infty
 \end{equation}
see Theorem \ref{mainth2} and Remark \ref{scal}.

Let us list our main assumptions
on the entropy function $F$ which will be used below.
A~typical example is given by
$F=\log$.

\begin{itemize}
\item [{ \bf A1)}]
$F$ is concave,increasing and $F(1)=0$
\item [{\bf A2)}]
$ \lim_{y \to 0} y F(y) =0$, $\lim_{y \to \infty} F(y)=\infty$
\end{itemize}
\begin{itemize}
\item [{\bf A3)}]
$y F(y)$ is convex on $[0,1+\Delta]$  for some $\Delta>0$
\item [{\bf A4)}]
there exists $y_0 \ge 1$ such that
$y F'(y) $ is non-increasing and $yF'(y) \le 1$ on $[y_0,\infty)$.
\end{itemize}

{\bf Remark.} {\it Assumptions A1) and A2) will be used throughout the paper. Assumptions A3) and A4) will be used for tight estimates. }

In Section 3 we obtain sufficient conditions for the related tight inequalities.
The case of the $F$-inequality follows immediately from the main result (Theorem \ref{Ftight}) without
any further assumptions. In the case of
modified log-Sobolev inequalities  we restrict ourselves to a
special choice of a cost function. Namely, we consider  for every
$1<\alpha  \le 2$  the corresponding family of cost functions
$c_{A,\alpha}$ given by (\ref{cost}).
Under some additional assumptions on the entropy, we prove a modification of (\ref{main}), where
$\int_{\R^d} f^2 \,d \mu$ is replaced by $\mbox{Var}_{\mu} f$ (Theorem \ref{modtight}). In the proof
we use techniques developed \cite{GGM}.

Before we give the precise formulation of the main result of Sections 3 and 4, let us briefly explain
the relationships between functions $F$, $c$, and $\mathcal{I}_{\mu}$ leading to
tight inequalities. We want to prove (\ref{12.11}).
It turns out that under assumptions A1)-A4) on  $\varphi$
every entropy function $F$  such that $F \sim A\varphi^\tau$, $\tau \le 1$ satisfies
\begin{equation}
\label{equv}
\Phi(x)\le F^{-1}(1+x) \sim \varphi^{-1}\Bigl(\Bigl[\frac{x+1}{A}\Bigr]^{1/\tau} \Bigr).
\end{equation}
Assume, in addition, that $\mathcal{I}_{\mu}(t) \ge k t \varphi^{1-\frac{1}{\alpha}}(t)$ for some $1 < \alpha \le 2$.
Now take a cost function $c$  such   that $c \sim B |x|^{q}$. We set
$$
q=\frac{\tau}{\tau -1 + \frac{1}{\alpha}}.
$$
Then
$$
 \Phi\Bigl(\delta c\Bigl[\frac{t F(\frac{1}{t})}{{\mathcal I}_{\mu}(t)} \Bigr] \Bigr)
 \le F^{-1}\Bigl(1+\varepsilon(\delta)F (1/t)\Bigr),
 $$
 where $\lim_{\delta \to 0} \varepsilon(\delta) =0$.
 Taking into account  property A4),
 one can easily show that $F^{-1}\bigl(1+\varepsilon F(1/t)\bigr) \le a t^{-p}$ for some $p<1$ and sufficiently small $\varepsilon$.
 Hence (\ref{12.11}) holds.

We  consider the generalized entropies defined by
$$f \to
\int_{\R^d} f F_{\tau}\Bigl(\frac{f}{\mu(f)}\Bigr) \,d \mu,
$$ where
$$
F_{\tau}(x) =  \left\{
\begin{array}{lcr}
 \varphi (x)  \quad       \mbox{if} \ 0 < x \le x_0 \\
 \frac{1}{\tau} \bigl( \varphi^{\tau}(x) -1 \bigr) +1 \quad \mbox{if} \ x \ge x_0, \\
\end{array}
\right.
$$
$\varphi$ satisfies A1)-A4), $\tau \le 1$ and $x_0$ is chosen in such a way that $\varphi(x_0)=1$.

Recall that $m_{f} = \inf \{ t: \mu(f > t) \le 1/2 \}$ is called
the median of $f$. Throughout the paper we assume that $\mu$ has
convex support.

\begin{theorem}
\label{convex}
Let $\varphi$ satisfy $A1)-A4)$ and let
 ${\mathcal I}_{\mu}$ satisfy
$$
{\mathcal I}_{\mu}(t) \ge k t \varphi \Bigl( \frac{1}{t}\Bigr)^{1-\frac{1}{\alpha}}
$$
for some $k>0$, $1< \alpha \le 2$ and $t \le 1/2$.
Then for every $2 \bigl(1-\frac{1}{\alpha} \bigr) \le  \tau \le 1$
 there exists
$C_{\tau} > 0$ depending on $\tau, \alpha, k, \lambda_2, \Delta$, such that for every smooth $f$ one has
$$
\int_{\R^d} f^2 F_{\tau} \Bigl( \frac{f^2}{\int_{\R^d} f^2 \,d \mu}\Bigr) \,d \mu
\le
C_{\tau} \int_{\R^d} f^2 c_{A,\frac{\alpha \tau }{\alpha-1}} \Bigl( \frac{|\nabla f|}{|f|}\Bigr) \,d \mu.
$$
In particular,
$$
\int_{\R^d} f^2 F_{2\bigl(1-\frac{1}{\alpha}\bigr)}\Bigl(\frac{f^2}{\int_{\R^d} f^2 \,d \mu }\Bigr) \,d \mu
\le
C_{2 \bigl(1-\frac{1}{\alpha} \bigr)} \int_{\R^d} | \nabla f |^2 \,d \mu,
$$
$$
\int_{\R^d} f^2 \varphi\Bigl(\frac{f^2}{\int_{\R^d} f^2 \,d \mu }\Bigr) \,d \mu \le C_1 \int_{\R^d}
f^2 c^{*}_{A,\alpha} \Bigl( \frac{|\nabla f|}{|f|}\Bigr) \,d \mu =
C_1 \int_{\R^d} f^2 c_{A,\frac{\alpha}{\alpha-1}}
\Bigl( \frac{|\nabla f|}{|f|}\Bigr) \,d \mu.
$$

In particular, the result holds if $\mu$ is convex and $g: \R^{+}
\to \R $ is increasing
 such that $\int_{\R^d} e^{g(r)} \,d \mu =1$
and for some $C>0$ one has
$$
\frac{g(r)}{\varphi^{1-\frac{1}{\alpha}}(e^{g(r)})} \ge C r.
$$
\end{theorem}

Obviously, if  $\mu$ is convex, $\varphi=\log$ and
\begin{equation}
\label{exp-pow-int}
\int_{\R^d} e^{\varepsilon |x|^{\alpha}} \,d \mu < \infty
\end{equation}
for some $\varepsilon>0$, we obtain
$$
\mbox{Ent}_{\mu} f^2 \le C_1 \int_{\R^d} f^2 c_{A,\frac{\alpha}{\alpha-1}}
\Bigl( \frac{|\nabla f|}{|f|}\Bigr) \,d \mu.
$$
In particular, we generalize Wang's criterion for convex
measures as well as  the  result of \cite{GGM}. Note that unlike
\cite{GGM} we deal directly with multidimensional distributions
and use a slightly different  cost function for $d \ge 2$.
We also apply the method developed in Theorem \ref{mainth}
to establish the following result (Theorem \ref{Bob-Zeg2}): let $\mu$ be a convex measure satisfying
(\ref{exp-pow-int}) for some $\alpha>1$. Then
$$
\mbox{\rm Ent}_{\mu} |f|^{\beta} \le C \Bigl[ \int_{\R^d} |\nabla f|^{\beta} \,d\mu
+ \mbox{\rm Var}_{\mu} |f|^{\frac{\beta}{2}}
\bigr].
$$
This inequality is  weaker than (\ref{Bob-Zeg}) but unlike (\ref{Bob-Zeg}) it is established for arbitrary convex measure.

During the preparation of the paper the author learned from Franck
Barthe that modified Sobolev inequalities for convex measures can
be obtained by using the transfer principle method (see \cite{Bar})
and the results from \cite{GGM}. However, this requires to prove
fist inequalities on the real line by different methods.  Another
achievement in this direction has been obtained
 by Nathael Golzan in
\cite{Goz2}, where he has proved a criterion for transportation
inequalities of Talagrand type for the real line. In particular,
his result implies  modified Sobolev inequalities for convex
measures on the real line, since they are known to be equivalent
to transportation inequalities in the log-concave case. The author
thanks the anonymous referee for very helpful comments.
\section{Main result}

Consider  a probability measure $\mu=\rho \,dx$ on $\R^d$.
We assume throughout that $X:=\mbox{supp} (\mu)$  is convex. In addition, without
loss of generality we
assume that $0 \in X$.
Set:
$$B_{r} = \{x: |x| \le r\}.
$$
We denote by $R(X) \in (0,\infty]$ the smallest number such that  $X \subset B_{R(X)}$.
Recall that for every measurable mapping $F:X \to Y$
the image measure
$\mu_F$ on $Y$ is defined by
$$
\mu_{F}(A) = \mu \bigl( \{ x: F(x) \in A \} \bigr)
$$
for every Borel set $A \subset Y$.
For every non-negative function $f$ we denote by $\tilde{f}$ the corresponding
spherical rearrangement, i.e., the function of the form
$\tilde{f}(x) = g(|x|)$ such that  $g$ is increasing and
$$
\mu \circ f^{-1} = \mu \circ \tilde{f}^{-1}.
$$
This can be rewritten as
$$
\mu_{f} = \mu_{r} \circ g^{-1}
$$
where $\mu_{f} = \mu \circ f^{-1}$ and $\mu_{r}$ is the image of $\mu$  under $x \to |x|$.
For a probability measure $\nu$ on $\R^{+}$ let us set
$$
F_{\nu}(t) = \nu([0,t))
$$
and
$$
G_{\nu}(u) = \{ \inf s: F_{\nu}(s) \ge u \}.
$$
Then $g$ has the form
\begin{equation}
\label{repr}
g = G_{\mu_{f}} \circ F_{\mu_r}.
\end{equation}
We denote by $B^{c}_r$ the complement of $B_r$ and
by $R_t>0$
the number such that
$$
\mu(|x| \le R_t)=t, \ R_1 = R(X).
$$
Since $X$ is convex and $0 \in X$, $R_t$ is well-defined.

For every $F:\R^{+} \to \R$ we define the corresponding isoperimetric  function $I_F$.
First we set
$$
J_{F}(s) = \frac{s F \bigl(\frac{1}{s}\bigr)}
 {\mathcal{I}_{\mu}(s)} .
$$
Equivalently,
$$
J_{F}(s) = \sup_{A: A \subset \R^d, \ \mu(A)=s }
\Bigl[ \frac{s F \bigl(\frac{1}{s}\bigr)}
 {\mu^{+}(A)} \Bigr].
$$
Then we define
$$
I_F(r) = J_F (1-\mu(B_r) ).
$$
This is equivalent to
$$
I_{F}(r) = \sup_{A \in \MM_r} \frac{\mu(A) F \bigl( \frac{1}{\mu(A)} \bigr) }{\mu^{+}(A)},
$$
where $\MM_r = \{A : \mu(A) = 1-\mu(B_r)\}$.
We follow the agreement that $I_F(R(X)) = 0$.

In what follows we consider a convex cost function $c : \R^{+} \to \R^{+}$.
Let
$$
c^{*}(x) = \sup_{y \in \R^{+}} \bigl( \bigl<x,y\bigr> - c(y) \bigr).
$$
We recall that $c$ is called superlinear if $\lim_{x \to \infty} \frac{c(|x|)}{|x|} = \infty$.
In what follows,   for simplicity we set
  $\mu(f^2) = \int_{\R^d} f^2 \,d \mu$.

\begin{theorem}
\label{mainth}
Let $c: \R^{+} \to \R^{+}$ be a convex  superlinear function such that $c(0)=0$
and let $F$ be  a function on $\R^{+}$ satisfying assumptions A1) and A2).
Let $K>1$. Assume that for $R=R_{\frac{K-1}{K}}$
one has
\begin{equation}
\label{mainassum}
  \int_{B^{c}_{R}} \Phi  \Bigl(4 c \circ I_F(|x|)\Bigr) \,d\mu < \infty,
\end{equation}
where
$$
\Phi(x) = \sup_{ y \in \R^{+}} \bigl( \bigl<x,y\bigr> - y F(y) + y \bigr)
= \bigl( yF(y) - y \bigr)^{*}(x).
$$
Then  there exist $B>0$, $C>0$  such that for every smooth
$f$ the following estimates hold:
\begin{equation}
\label{defmodLSI}
\int_{\R^d} f^2 F\Bigl(\frac{f^2}{\int_{\R^d} f^2 \,d \mu}\Bigr) \,d \mu
\le
4 \int_{\{f^2 \ge K\int f^2 \,d \mu\}} f^2 c^{*} \Bigl( \frac{|\nabla f|}{|f|}\Bigr) \,d \mu
+
B  \int_{\R^d} f^2 \,d\mu,
\end{equation}

\begin{align}
\label{defmodLSI2}
 \int_{\R^d} f^2 F\Bigl(\frac{f^2}{\int_{\R^d} f^2 \,d \mu}\Bigr) \,d \mu &
\\&
\nonumber
\le
C \int_{\R^d}
\bigl(f-\mu(f)\bigr)^2 c^{*} \Bigl( \frac{|\nabla f|}{\bigl|f - \mu(f)\bigr|}\Bigr) \,d \mu
+
B \cdot\mbox{\rm Var}_{\mu} f.
\end{align}

\end{theorem}

 \begin{proof}
  Let us fix some a Lipschitz function $f$. Without loss of generality we
  may assume that $f \ge \varepsilon>0$.
 Set  $\nu: = g \cdot   \mu$, where $g = F \bigl(\frac{f^2}{\int f^2 \,d \mu}\bigr)$.
 By a well-known result from measure theory one has
 \begin{align*}
 &
\int_{\R^d} f^2 F\Bigl(\frac{f^2}{\int_{\R^d} f^2 \,d \mu}\Bigr) \,d \mu =
 \int_{\R^d} f^2 g \,d \mu
 =
 \int_{\R^d} f^2  \,d \nu
 =
 \int_{0}^{\infty} \nu( f^2(x) > t) \,dt
 \\&
 =
 \int_{0}^{\infty} \Bigl( \int_{\{x: f^2(x) > t\}} g \,d \mu \Bigr) \,dt.
 \end{align*}
 We split this integral in the following two parts:
 $$
 I_1 =   \int_{0}^{K \mu(f^2)} \Bigl( \int_{\{x: f^2(x) > t \}} g \,d \mu \Bigr) \,dt,
 \quad
 I_2 =   \int_{K \mu(f^2) }^{\infty} \Bigl( \int_{\{x: f^2(x) > t \}} g \,d \mu \Bigr) \,dt.
 $$

 The following proof will be divided in several steps.

 {\bf Step 1}.
 Estimation of $I_1$.
 We show that for some $C(K)>0$ one has
 $$
 I_1 \le C(K) \mbox{Var}_{\mu} f.
 $$
 This part is quite elementary.
 By the concavity of $F$ one has
 $$
 g \le F'(1) \Bigl(\frac{f^2}{\mu(f^2)} -1\Bigr).
 $$
 Hence
 \begin{align*}
 \frac{I_1}{F'(1)} & \le \frac{1}{\mu(f^2)} \int_{\R^d}
 \Bigl( f^2 - \mu(f^2) \Bigr)
 \Bigl( \int_{0}^{K \mu(f^2)} I_{\{x: f^2(x) > t \}} \,dt \Bigr) \,d \mu.
 \\&
 = \frac{1}{\mu(f^2)} \int_{\R^d}
 \Bigl( f^2 - \mu(f^2) \Bigr)
 \min \Bigl( f^2, K \mu(f^2) \Bigr) \,d \mu
 \\&
 =  \frac{1}{\mu(f^2)} \int_{\R^d}
 \Bigl( f^2 - \mu(f^2) \Bigr)
 \Bigl[ \min \Bigl( f^2, K \mu(f^2) \Bigr) -   \mu(f^2)  \Bigr]\,d \mu.
 \end{align*}
 The latter equals
 \begin{align*}
  \frac{1}{\mu(f^2)} \int_{\{f^2 \le K \mu(f^2)\}}
 \Bigl( f^2 - \mu(f^2) \Bigr)^2  \,d \mu
 +
 (K-1)  \int_{\{f^2 \ge K \mu(f^2)\}}
 \Bigl( f^2 - \mu(f^2) \Bigr) \,d \mu.
 \end{align*}
 The first term can be estimated in the following way:
 \begin{align*}
   \frac{1}{\mu(f^2)} \int_{\{f^2 \le K \mu(f^2)\}} &
 \Bigl( f^2 - \mu(f^2) \Bigr)^2  \,d \mu
 \\&
 \le
    \frac{2}{\mu(f^2)} \int_{\{f^2 \le K \mu(f^2)\}}
 \Bigl( f^2 - \mu(f)^2 \Bigr)^2  \,d \mu
 +
 \frac{2}{\mu(f^2)} \Bigl[ \mbox{\rm{Var}}_{\mu} f \Bigr]^2
 \\&
 \le
     4(K+1)^2 \int_{\R^d}
 \Bigl( f - \mu(f) \Bigr)^2  \,d \mu
 +
 2 \mbox{\rm{Var}}_{\mu} f
 =
 (4(K+1)^2+2) \mbox{\rm{Var}}_{\mu} f.
 \end{align*}
 Further we get
 $$
  \int_{\{f^2 \ge K \mu(f^2)\}}
 \Bigl( f^2 - \mu(f^2) \Bigr) \,d \mu
 \le
   \int_{\{f^2 \ge K \mu(f^2)\}}
 \Bigl( f^2 - \mu(f)^2 \Bigr) \,d \mu .
 $$
 One can easily check that
 $$
 |f + \mu(f)| \le \frac{\sqrt{K}+1}{\sqrt{K}-1}  |f - \mu(f)|
 $$
 on $\{f^2 \ge K \mu(f^2)\}$.
 Hence
 \begin{equation}
 \label{tail-var}
 \int_{\{f^2 \ge K \mu(f^2)\}} \Bigl( f^2 - \mu(f^2) \Bigr) \,d \mu
 \le
 \frac{\sqrt{K}+1}{\sqrt{K}-1} \mbox{Var}_{\mu} f.
 \end{equation}
 Finally we obtain
 $$
 I_1 \le \Bigl[ (4(K+1)^2+2) + (\sqrt{K}+1)^2 \Bigr] F'(1) \mbox{\rm{Var}}_{\mu} f.
 $$

 {\bf Step 2}.
 Here we estimate   $I_2$ by a quantity depending on the isoperimetric function $I_{F}$.
 Let us set
 $$
 A_t =\{x: f^2(x) > t\}.
 $$
By the concavity of $F$ one has
\begin{align*}
I_2 & =
  \int_{K \mu(f^2)}^{\infty}
  \int_{\R^d} I_{A_t} F \Bigl( \frac{f^2}{\int_{\R^d} f^2 \,d \mu} \Bigr) \,d \mu
\,dt
 \\& \le
  \int_{K \mu(f^2)}^{\infty}  \mu(A_t ) \Bigl[  F \Bigl( \int_{A_t}\frac{f^2}{ \mu(A_t)
  \int_{\R^d} f^2 \,d \mu} \Bigr)\,d \mu \Bigr]
 \,dt
 \\&
 \le
 \int_{K \mu(f^2)}^{\infty}  \mu(A_t) F \Bigl( \frac{1}{ \mu(A_t )}\Bigr) \,dt
 \\&
 =
 \int_{K \mu(f^2)}^{\infty}  \mu(\{x: f^2(x) > t\} ) F \Bigl( \frac{1}{ \mu(\{x: f^2(x) > t\} )} \Bigr)
 \,dt.
 \end{align*}
Since $f$ is continuous and $X$ is convex,  the function
$t\mapsto \mu(A_t)$ is strictly decreasing on
$$[ \inf_{x \in X} f^2(x), \sup_{x \in X} f^2(x)].$$
 Hence one can find  a nondecreasing function $r_{f^2}(s)$ such that
 $$
 \mu(A_s) = \mu(B^c_{r_{f^2}(s)})
 $$
 and $r_{f^2}(0)=0$, $r_{f^2}(s)= R(X)$, if $s \ge \sup f^2$.
  Set
 $$
 f_h(x) = \sup_{\{|x-y| \le h\}} f(y).
 $$
 By the definition of $I_{F}$ we have
 \begin{align*}
 &
 I_2 \le
 \int_{K \mu(f^2)}^{\infty} I_{F} \bigl(r_{f^2}(t)\bigr) \mu^{+}(A_t) \,dt
 \\&
 \nonumber
 \le
 \underline{\lim}_{h \to 0+} \int_{K \mu(f^2)}^{\infty} I_{F}(r_{f^2}(t))  \frac{\mu(A^h_t) - \mu(A_t)}{h} \,dt,
 \end{align*}
 where $\{ x \in \R^d: f^2_h(x) > t\} = \{x \in \R^d: f^2(x)>t\}^{h}  = A^{h}_{t} $.
  Assume for a while that $s \to I_F(r_{f^2}(s))$ is locally integrable and define
 $$
 Z(t) : = \left\{
\begin{array}{lcr} \int_{K \mu(f^2)}^{t} I_{F}(r_{f^2}(s)) \,ds, \ t  \ge K \mu(f^2)
\\ 0, \  t  \le K \mu(f^2).
\end{array}
\right.
 $$
Applying the formula
$$
\int \Phi(f^2) \,d \mu = \int_{0}^{\infty} \Phi'(t) \mu(A_t) \,d t ,
$$
which holds for every increasing $\Phi$ such that $\Phi(0)=0$, we get
\begin{align*}
I_2
&\le \underline{\lim}_{h \to 0+} \int_{\R^d} \frac{Z(f^2_h) - Z(f^2)}{h} \,d \mu
\\&
\le
2\int_{\{f^2 \ge K \mu(f^2)\}} I_{F}\bigl(r_{f^2}(f^2)\bigr) |f| |\nabla f| \,d \mu.
\end{align*}
It remains to note that this estimate still holds even if $I_{F}(r_{f^2})$ is not locally integrable.
Indeed, approximating $I_{F}$ by $I^{N}_{F} = I_{F} \wedge N$,
we obtain in the same way as above that
\begin{align*}
\int_{K \mu(f^2)}^{\infty} I_{F}^ N \bigl(r_{f^2}(t)\bigr) \mu^{+}(A_t) \,dt
& \le
 2 \int_{\{f^2 \ge K \mu(f^2)\}} I_{F}^N  \bigl(r_{f^2}(f^2)\bigr) |f| |\nabla f| \,d \mu
 \\&
 \le
 2 \int_{\{f^2 \ge K \mu(f^2)\}} I_{F}  \bigl(r_{f^2}(f^2)\bigr) |f| |\nabla f| \,d \mu.
\end{align*}
We apply the monotone convergence theorem
$$
I_2 \le  \int_{K \mu(f^2)}^{\infty} I_{F} \bigl(r_{f^2}(t)\bigr) \mu^{+}(A_t) \,dt
= {\lim}_N \int_{K \mu(f^2)}^{\infty} I_{F}^ N \bigl(r_{f^2}(t)\bigr) \mu^{+}(A_t) \,dt,
$$
and obtain the claim.

{\bf Step 3.}
Estimation of
$$
\int_{\{f^2 \ge K \mu(f^2)\}} I_{F}\bigl(r_{f^2}(f^2)\bigr) |f| |\nabla f| \,d \mu.
$$
We complete the desired estimate by using the Young inequality. In this part rearrangement
techniques will be employed.
Namely, in the estimate below  we replace
$I_{F}\bigl(r_{f^2}(f^2)\bigr)$
by $I_{F}\bigl(r_{f^2}(\tilde{f}^2)\bigr)$  and take into account that $r_{f^2}(\tilde{f}^2(x))=|x|$
on the set $\{ x: |\nabla f(x)| \ne 0 \}$.

Let $\R_{\delta} = \{ t: \mu \circ \bigl(f^2\bigr)^{-1} (t) >0\}$ be the set of atoms of
the measure $\mu \circ \bigl( f^2\bigr)^{-1}$.
Note that $| \nabla f| =0$ almost everywhere on $D = \{x : f^2(x) \in \R_{\delta} \}$.
Hence by the Young inequality we find
\begin{align}
\label{young}
&
2\int_{\{f^2 \ge K \mu(f^2)\}} I_{F}(r_{f^2}(f^2)) |f| |\nabla f| \,d \mu
\le
2 \int_{\{f^2 \ge K \mu(f^2)\}} f^2 c^{*} \Bigl( \frac{|\nabla f|}{|f|}\Bigr)
\,d \mu
\\&
\nonumber
+ 2 \int_{\{f^2 \ge K \mu(f^2)\} \cap D^{c}}f^2  \bigl[ c \circ  I_F \bigl(r_{f^2}(f^2)\bigr) \bigr] \,d \mu.
\end{align}
 Let $O_K = \{x: f^2(x) \ge K \mu(f^2)\} \cap D^{c}$.
One has $$I_{O_K} = I_{\{f^2 \ge K \mu(f^2)\}} \cdot I_{\R^{c}_{\delta}} (f^2)$$
and by the Young inequality
\begin{align*}
&
2 \int_{O_{K}} f^2  c \bigl(  I_F(r_{f^2}(f^2)) \bigr) \,d \mu
 =
 2 \int_{\R^d} f^2  I_{O_{K}} c \bigl(  I_F(r_{f^2}(f^2)) \bigr) \,d \mu
\\&
=
 \frac{1}{2} \mu(f^2) \int_{\R^d}
 \Bigl[\frac{f^2}{ \mu(f^2)} \Bigr]
 \Bigl[4 I_{O_{K}} c \bigl( I_F(|r_{f^2}(f^2)|) \bigr)  \Bigr]\,d \mu
 \,d \mu
 \\&
 \le
\frac{1}{2} \int_{\R^d} f^2 \Bigl[ F \Bigl( \frac{f^2}{ \mu(f^2)}\Bigr) -1 \Bigr]\,d \mu
+ \frac{1}{2}\mu(f^2) \int_{\R^d}
  \Phi  \Bigl(4 I_{O_K} c \bigl( I_F(r_{f^2}(f^2))\bigr)  \Bigr)
  \,d \mu.
\end{align*}
Since $f$ and $\tilde{f}$ have the same laws considered as  random variables on the
probability space $(\R^d, \mu)$, one has
$$
\int_{\R^d}
  \Phi  \Bigl(4 I_{O_{K}} c \bigl( I_F(r_{f^2}(f^2))\bigr)   \Bigr)
  \,d \mu
  =
\int_{\R^d}
  \Phi  \Bigl(4 I_{\tilde{O}_{K}} c \bigl( I_F(r_{f^2}(\tilde{f}^2))\bigr)    \Bigr)
  \,d \mu
$$
where $\tilde{O}_{K} = \{x: \tilde{f}^2(x) \ge K \mu(f^2)\} \cap \{ x: \tilde{f}^2(x) \in \R^{c}_{\delta}\}$.
By the definition of  $\tilde{f}$ we have
$$
 \mu\bigl(\{ y: f^2(y) > \tilde{f}^2(x) \}\bigr)
= \mu\bigl(\{ y: {\tilde f}^2(y) > \tilde{f}^2(x) \}\bigr).
 $$
Then for every such $x$ by the definition of $r_{f^2}$
we have
 $$
 \mu \Bigl(B^c_{r_{f^2}({\tilde f}^2)}\Bigr) =\mu\bigl(\{ y: {\tilde f}^2(y) > \tilde{f}^2(x) \}\bigr)
 = \mu(\{ y: |y| > |x|\}).
 $$
 Indeed, otherwise there exist $r_1 < r_2$ such that $\tilde{f}(z) = \tilde{f}(x)$
for every $z: r_1 \le |z| \le r_2$. But this implies that
$\mu \bigl(y: f(y) = \tilde{f}(x) \bigr)>0$.
 Hence $r_{f^2}({\tilde f}^2)(x) = |x|$ on $\tilde{O}_{K}$.
Moreover, if $x \in \{\tilde{f}^2 \ge K \mu(\tilde{f}^2) \}$, then by the Chebyshev inequality
$$
 \mu\bigl(B^c_{|x|}\bigr) = \mu \Bigl(B^c_{r_{f^2}({\tilde f}^2)}\Bigr)
 \le  \mu \bigl( \{\tilde{f}^2 \ge K \mu(\tilde{f}^2) \}\bigr)
 \le 1/K.
 $$
 Hence
$|x| = r_{f^2}({\tilde f}^2(x)) \ge R_{(K- 1)/K}$ if $ x \in \{\tilde{f}^2 \ge K \mu(\tilde{f}^2) \}$.
Thus
$$\tilde{O}_{K} \subset \{ x: |x| \ge R_{(K- 1)/K} \}.
$$
Hence
\begin{align*}
\int_{\R^d}
  \Phi  \Bigl(4 I_{\tilde{O}_{K}} c \bigl( I_F(r_{f^2}(\tilde{f}^2))\bigr)    \Bigr)
  \,d \mu
  \le
 \Phi(0) + \int_{B^c_{R_{(K- 1)/K}}}
  \Phi  \Bigl(4  c \bigl( I_F(|x|)\bigr)    \Bigr)
  \,d \mu := \tilde{B} < \infty.
\end{align*}
Finally
\begin{align*}
\frac{1}{2} \int_{\R^d} f^2 \Bigl[ F \Bigl( \frac{f^2}{ \mu(f^2)}\Bigr) -1 \Bigr]\,d \mu
+
& \frac{\mu(f^2)}{2}\int_{\R^d}
  \Phi  \Bigl(4 I_{O_{K}} c \bigl( I_F(r_{f^2}(f^2))\bigr)   \Bigr)
  \,d \mu
 \\& \le
 \frac{1}{2} \int_{\R^d} f^2  F \Bigl( \frac{f^2}{ \mu(f^2)}\Bigr) \,d \mu
 +
 \frac{\tilde{B} -1}{2} \int_{\R^d} f^2 \,d \mu.
\end{align*}
and
$$
I_2  \le \frac{\tilde{B} -1}{2}  \int_{\R^d} f^2 \,d \mu
+
\frac{1}{2} \int_{\R^d} f^2  F \Bigl( \frac{f^2}{ \mu(f^2)}\Bigr)
\,d \mu
+
2 \int_{\{f^2 \ge K \mu(f^2)\}} f^2 c^{*} \Bigl( \frac{|\nabla f|}{|f|}\Bigr)
\,d \mu
.
$$
Combining all the inequalities obtained above, we get (\ref{defmodLSI}).

 The proof of (\ref{defmodLSI2}) is similar and we just briefly describe the main difference.
 Instead of (\ref{young}) we use
 \begin{align*}
&
2\int_{\{f^2 \ge K \mu(f^2)\}} I_{F}(r_{f^2}(f^2)) |f| |\nabla f| \,d \mu
\le
\\&
C' \int_{\{f^2 \ge K \mu(f^2)\}} \bigl(f-\mu(f)\bigr)^2
c^{*} \Bigl( \frac{|\nabla f|}{|f-\mu(f)|}\Bigr)
\,d \mu
\\&
\nonumber
+ C' \int_{O_{K}} \bigl(f-\mu(f)\bigr)^2   \bigl[ c \circ  I_F \bigl(r_{f^2}(f^2)\bigr) \bigr] \,d \mu.
\end{align*}
This follows from the Young inequality and the observation that
 $$
 f^2 \le \frac{K}{(\sqrt{K}-1)^2}  \bigl( f - \mu(f)\bigr)^2
 $$
 on $\{f^2 \ge K \mu(f^2)\}$.
 In the same way as above we estimate the second term by
 $\mbox{\rm Var}_{\mu} f$
 and  $\int_{\R^d} {\tilde f}^2 F\Bigl(\frac{\tilde{f}^2}{\int_{\R^d} \tilde{f}^2 \,d \mu}\Bigr) \,d \mu$,
 where $\tilde{f} = f -\mu(f)$.
 Finally, by (\ref{defmodLSI})
one has
 $$\int_{\R^d} {\tilde f}^2 F\Bigl(\frac{\tilde{f}^2}{\int_{\R^d} \tilde{f}^2 \,d \mu}\Bigr) \,d \mu
 \le
4 \int_{\R^d} (f-\mu(f))^2 c^{*} \Bigl( \frac{|\nabla f|}{|f-\mu(f)|}\Bigr) \,d \mu
+
B  \cdot \mbox{Var}_{\mu} f.
 $$
The proof is complete.
\end{proof}

\begin{example}
Assume that $c$ is a convex superlinear function satisfying (H). Let $\mu$ be a convex measure such that
$\int_{\R^d} e^{\varepsilon c(r)} \,d \mu < \infty$ for some $\varepsilon >0$.
Then for every $K$ there exist $B, C >0$ such that
\begin{equation}
\mbox{{\rm Ent}}_{\mu} f^2
\le
C \int_{\{f^2 \ge K \int f^2 \,d \mu\}} f^2 c^{*} \Bigl( \frac{|\nabla f|}{|f|}\Bigr) \,d \mu
+
B \int_{\R^d} f^2 \,d\mu.
\end{equation}
\end{example}
\begin{proof}
Let $F = \log$.
 It will be shown below that
$\sup_{r \ge R_{1/2}}\frac{I_{\log}(r)}{r} < \infty$ for every
convex $\mu$ (Lemma \ref{conv-iso}). The result follows
immediately from Theorem \ref{mainth}.
\end{proof}

\begin{theorem}
\label{mainth2}
Let $c: \R^{+} \to \R^{+}$ be a convex  superlinear function such that $c(0)=0$.
Assume that $F$ satisfies assumptions A1)-A2)
and there exists $K>1$ such that
\begin{equation}
\label{1dimint}
\int_{0}^{1/K} \Phi \Bigl( 4c  \Bigl[ \frac{t F(\frac{1}{t})}{{\mathcal I}_{\mu}(t)} \Bigr] \Bigr) \,dt < \infty.
\end{equation}
Then inequalities (\ref{defmodLSI}) and (\ref{defmodLSI2}) hold.
\end{theorem}
\begin{proof}
By the definition $I_{F}$ one has
$$
I_{F}(r) = \frac{(1-\mu(B_r)) F \bigr(\frac{1}{1-\mu(B_r)} \bigr)}{{\mathcal I}_{\mu}(1-\mu(B_r))}.
$$
It suffices to show that
$$
  \int_{B^{c}_{R_{\frac{K-1}{K}}}} \Phi  \Bigl(4 c \circ I_F(|x|)\Bigr) \,d\mu < \infty.
$$
The mapping $ \R^d \ni x \to 1 - \mu(y: |y| \le |x|)=t \in [0,1]$
transforms $\mu$ into Lebesgue
measure on  $[0,1]$. Hence the integrability of $\Phi(4 c(I_F))$ is equivalent to
(\ref{1dimint}) for some $\varepsilon >0$.
\end{proof}

\begin{remark}
\label{scal}
Note that the constant $4$ in (\ref{mainassum}) and (\ref{1dimint}) provides
 yields the term $$4  \int_{\{f^2 \ge K \int f^2 \,d \mu\}} f^2 c^{*}
  \Bigl( \frac{|\nabla f|}{|f|}\Bigr) \,d \mu$$ in (\ref{defmodLSI}). However, if $c$ satisfies (H),
it is more convenient to assume that
$$
\int_{0}^{1/K} \Phi \Bigl( \delta c  \Bigl[ \frac{t F(\frac{1}{t})}{{\mathcal I}_{\mu}(t)} \Bigr] \Bigr) \,dt < \infty
$$
for some $\delta>0$, $K>1$. It is easy to check (just apply Theorems \ref{mainth}, \ref{mainth2}
to $\tilde{c} = \varepsilon c$ with appropriate $\varepsilon$) that (\ref{defmodLSI},
(\ref{defmodLSI2}) still hold (eventually with some
other constant in place of~$4$).
\end{remark}

The following theorem is a direct corollary of (\ref{defmodLSI2}).

\begin{theorem}
\label{Ftight}
Let  $F$ and $\mu$ satisfy the assumptions of
Theorem \ref{mainth} with  $c=\delta |x|^2$
and some $\delta>0$.
Then for every smooth $f$ one has
$$
\int_{\R^d} f^2 F\Bigl(\frac{f^2}{\int_{\R^d} f^2 \,d \mu}\Bigr) \,d \mu
\le
C \int_{\R^d} |\nabla f|^2 \,d \mu
+
B \cdot\mbox{\rm Var}_{\mu} f.
$$
In particular, the result holds if  assumptions A1)-A2) are fulfilled
and there exist $K>1$, $\delta >0 $ such that
\begin{equation}
\int_{0}^{1/K} \Phi \Bigl( \delta  \Bigl[ \frac{t F(\frac{1}{t})}{{\mathcal I}_{\mu}(t)} \Bigr]^2 \Bigr) \,dt < \infty.
\end{equation}
\end{theorem}

\begin{example} (d=1)
Consider a probability measure on the real line $\mu = e^{-V(t)}
\,dt$. In the one-dimensional case the proof can be simplified. We
omit here the detailed proof and just briefly explain the main
ideas. Instead of using the coarea inequality one can
apply the
Newton--Leibnitz formula
$$
f(x) = f(m) + \int_{m}^{x} f'(s) \,ds,
$$
where $m \in \R$. It is convenient to take for $m$ the median of $\mu$.
The use of the Newton--Leibnitz formula allows to
use the  simplified analog of the isoperimetric function ${\tilde{\mathcal{I}}_{\mu}}$.
Let $ 0 \le t \le 1/2$. Define $u(t) \le m$ and $v(t) \ge m$ as follows:
$$
\mu((-\infty, u(t)]) =  \mu([v(t),\infty)=t.
$$
Then
$$
{\tilde{\mathcal{I}}_{\mu}} (t) =
\min \bigl\{
e^{V(u(t))},   e^{V(v(t))}\bigr\} t.
$$
One can get the following analog of Theorem \ref{Ftight}:

{ \rm
Let  assumptions A1)-A2) be satisfied
and let $K>2$ and $\delta >0$ be such that
\begin{equation}
\int_{0}^{1/K} \Phi \Bigl( \delta  \Bigl[ \frac{t F(\frac{1}{t})}
{\tilde{{\mathcal I}}_{\mu}(t)} \Bigr]^2 \Bigr) \,dt < \infty.
\end{equation}
 Then
\begin{equation}
\label{05.11.06}
\int_{\R} f^2 F\Bigl(\frac{f^2}{\int_{\R} f^2 \,d \mu}\Bigr) \,d \mu
\le
C \int_{\R} |f'|^2 \,d \mu
+
B \cdot\int_{\R} (f - f(m))^2 \,d\mu
\end{equation}
for some $B, C>0$ and every smooth $f$.

If, in addition, $\mu$ satisfies the Poincar{\'e} inequality,  the term
$\int_{\R} (f - f(m))^2 \,d\mu$ can be estimated by $C' \int_{\R} |f'|^2 \,d \mu$
(see \cite{BZ}) and be omitted in (\ref{05.11.06}):
\begin{equation}
\label{05.11.06(2)}
\int_{\R} f^2 F\Bigl(\frac{f^2}{\int_{\R} f^2 \,d \mu}\Bigr) \,d \mu
\le
C \int_{\R} |f'|^2 \,d \mu.
\end{equation}
}

As an example consider the following measure on
the line:
$$
\mu= Z e^{-|x| \log(1+x^2)} \,dx.
$$
It can be easily verified that as $s \to \infty$ one has
$$
\mu((-\infty,-s]) = \mu([s,\infty)) \sim  \frac{Ze^{-|s| \log(1+s^2)}}{\log(1+s^2)}.
$$
Since $\mu^{+}([s,\infty)) = Ze^{-|s| \log(1+s^2)}$, we get
$$
{\tilde{{\mathcal I}}_{\mu}(t)}
\ge C' t \log\bigl (\log(1/t)) \bigr)
$$
for some $C'$ and every $t\ge 1/2$.
Let us choose a function $F$ satisfying assumptions A1)-A2) of Theorem \ref{mainth} such that
$$
F(x) \sim \log^2(\log x).
$$
for large values of $x$. In this case
$$
\Phi(y) \sim \exp(e^{\sqrt{y}})
$$
for large $y$.
Hence for any sufficiently small $\delta$ and all
$t\in [0,1/2]$ one has
$$\Phi \Bigl( \delta  \Bigl[ \frac{t F(\frac{1}{t})}
{\tilde{{\mathcal I}}_{\mu}(t)} \Bigr]^2 \Bigr)
\le
\exp( \log^{p}(1/t) ),
$$
where $p$ can be done arbitrary small. Since
$$\int_{0}^{1/2} \exp( \log^{p}(1/t) )\,dt < \infty
$$
for $p <1$, we obtain (\ref{05.11.06(2)}).
\end{example}

\section{Tight estimates}
In this section we establish some tight estimates, i.e.,
 estimates
whose right-hand sides  vanish on constant functions.
The case of the $F$-Sobolev inequality has been already considered in
Theorem \ref{Ftight}.
Unlike the $F$-Sobolev inequality, the case of tight modified log-Sobolev inequalities
is more difficult. We use an idea from \cite{GGM} and consider two cases: the case
of large entropy and the case of small entropy.
The large entropy case follows immediately from our main result.
In the case of small entropy we reduce   the problem to the $F$-inequality.

In what follows we assume that there exists $\lambda_2 >0$ such that for every
smooth $f$ one has
\begin{equation}
\label{P2}
\int_{\R^d} \bigl( f - m_{f} \bigr)^2 \,d \mu \le \lambda_2 \int_{\R^d} |\nabla f|^2 \,d \mu.
\end{equation}

Since $\int_{\R^d} \bigl( f - \int_{\R^d} f \,d \mu \bigr)^2 \,d \mu \le \int_{\R^d}
\bigl( f - m_{f} \bigr)^2 \,d \mu$, this inequality is stronger than the classical $L^2$-Poincar{\'e} inequality.

\begin{definition}
We say that a probability measure $\mu$ satisfies the
Cheeger isoperimetric inequality if
there exists $\lambda_1 >0$ such that for every Borel set
$A$ one has
\begin{equation}
\label{Cheeger}
\min(\mu(A),1-\mu(A)) \le \lambda_1 \mu^{+}(A).
\end{equation}
\end{definition}

Inequality (\ref{Cheeger})  is equivalent to the following $L^1$-Poincar{\'e}-type inequality:
\begin{equation}
\label{l1poin}
\int_{\R^d} \bigl|f - \int_{\R^d} f \,d \mu\bigr| \,d \mu \le \lambda_1 \int_{\R^d} | \nabla f| \,d \mu.
\end{equation}
It was shown in \cite{BH} that  (\ref{Cheeger}) implies (\ref{P2}).
It is known that every convex measure satisfies (\ref{l1poin})
with some $\lambda_1$ (see \cite{KLS} and \cite{BH}).

We start this section with several lemmas.

\begin{lemma}
\label{dual-estimates}
Let $F$ satisfy assumptions A1), A2) and A4)
Then for every $\delta\in (0,1/2]$,
 there exists  $T$ depending on $\delta$ and $y_0$
 such that for any $y \ge T$ one has
$$
\Phi\Bigl( \delta F\bigl(y) \Bigr)\le y^{2\delta}.
$$
\end{lemma}
\begin{proof}
Since $F$ is increasing and $\lim_{y \to \infty} F(y) = \infty$, the supremum of
$$
xy - y F(y) +y
$$
is attained at some $y^*$. Moreover, there exists $x_0$ such that
$y^* \ge y_0$ if $x \ge x_0$. In this case one has
\begin{equation}
\label{extr}
x = F(y^*) + y^* F'(y^*) -1
\end{equation}
and by the properties of $F$
$$
F(y^*) -1 \le x \le F(y^*).
$$
Consequently,
$$
y^{*} \le F^{-1}(1+x)
$$
and by (\ref{extr}) we find
$$
\Phi(x) = x y^{*} - y^* F(y^*) +y^* = (y^*)^2 F'(y^*).
$$
Hence for  any $x \ge x_0$ one has
\begin{equation}
\label{dual-est}
\Phi(x) \le y^{*} \le F^{-1}(1+x).
\end{equation}
Next, for any $y \ge y_0$, we have
$$
F(y^{2\delta}) - F(y^{2\delta}_0)
=
2\delta \int_{y_0}^{y} s^{2\delta-1}F'(s^{2\delta}) \,ds.
$$
Taking into account that $s^{2\delta} \le s$, we get by A4)
$$
F(y^{2\delta}) - F(y^{2\delta}_0) = 2\delta \int_{y_0}^{y}\frac{s^{2\delta}F'(s^{2\delta})}{s}\,ds
\ge
2\delta \int_{y_0}^{y}\frac{sF'(s)}{s}\,ds
=
2\delta \bigl( F(y) - F(y_0) \bigr).
$$
Finally,
$$
\delta F(y) \le \delta F(y_0)- \frac{1}{2}F(y^{2\delta}_0) + \frac{1}{2}F(y^{2\delta}).
$$
Thus, if $F(y) \ge \frac{x_0}{\delta}$, we obtain
by (\ref{dual-est}) that
$$
\Phi(\delta F(y))
\le
F^{-1} \Bigl(1+\delta F(y_0)- \frac{1}{2}F(y^{2\delta}_0) + \frac{1}{2}F(y^{2\delta})\Bigr).
$$
Choosing $T \ge F^{-1} \bigl(\frac{x_0}{\delta}\bigr)$ in such a
way that $\frac{1}{2}F(y^{2\delta}) \ge 1+\delta F(y_0)-
\frac{1}{2}F(y^{2\delta}_0)$ for $y \ge T$, we obtain
$$\Phi(\delta F(y))
\le
F^{-1} \Bigl(1+\delta F(y_0)- \frac{1}{2}F(y^{2\delta}_0) + \frac{1}{2}F(y^{2\delta})\Bigr)
\le
F^{-1} \bigl( F(y^{2\delta}) \bigr) \le y^{2\delta}.
$$
The proof is complete.
\end{proof}

\begin{lemma}
Let $\mu$ be a probability measure and let
$F$  satisfy  assumptions A1), A2), and A4).
Then there exists $C>0$ such that for all
$f,g \in L^2(\mu)$ one has
\label{entr-bounds}
$$
\int_{\R^d} f^2 F \Bigl( \frac{g^2}{\int_{\R^d} g^2 \,d \mu}\Bigr) \,d \mu
\le
2\int_{\R^d} f^2 F \Bigl( \frac{f^2}{\int_{\R^d} f^2 \,d \mu}\Bigr) \,d \mu +
C \int_{\R^d} f^2 \,d \mu.
$$
\end{lemma}
\begin{proof}
Set
$$
u = F \Bigl( \frac{g^2}{\mu(g^2)}\Bigr) +
\frac{g^2}{\mu(g^2)} F'\Bigl( \frac{g^2}{\mu(g^2)} \Bigr) -1
$$
Since $F'>0$, one has
\begin{align*}
\int_{\R^d} f^2 F \Bigl( \frac{g^2}{\mu(g^2)}\Bigr) \,d \mu
\le
\int_{\R^d} f^2 u \,d \mu + \int_{\R^d} f^2 \,d \mu.
\end{align*}
By the Young inequality
\begin{align*}
\int_{\R^d} f^2 u & \,d \mu
=
\int_{\R^d} f^2 \,d \mu \int_{\R^d} \frac{f^2}{\mu(f^2)} u \,d \mu
\\& \le
2\int_{\R^d} f^2 \,d \mu \Bigl( \int_{\R^d}  \Bigl[\frac{f^2}{\mu(f^2)} F \Bigl( \frac{f^2}{
\mu(f^2)}\Bigr)  - \frac{f^2}{\mu(f^2) } \Bigr]\,d \mu \Bigr)
+ 2\int_{\R^d} f^2 \,d \mu \int_{\R^d} \Phi(u/2) \,d \mu.
\end{align*}
Hence
\begin{align*}
&
\int_{\R^d} f^2 F \Bigl( \frac{g^2}{\mu(g^2)}\Bigr) \,d \mu
\le
2\int_{\R^d} f^2 \,d \mu \Bigl( \int_{\R^d}  \Bigl[\frac{f^2}{\mu(f^2)} F \Bigl( \frac{f^2}{
\mu(f^2)}\Bigr)  - \frac{f^2}{\mu(f^2) } \Bigr]\,d \mu \Bigr)
\\&
+ 2\int_{\R^d} f^2 \,d \mu \int_{\R^d} \Phi(u/2) \,d \mu + \int_{\R^d} f^2 \,d \mu
\\&
= 2\int_{\R^d} f^2 F \Bigl( \frac{f^2}{ \mu(f^2)}\Bigr) \,d \mu
+
\int_{\R^d} f^2 \,d \mu \Bigl(\int_{\R^d}(2\Phi(u/2)-1) \,d \mu \Bigr).
\end{align*}
Using estimate $\Phi(x) \le F^{-1}(1+x)$  obtained in the proof of Lemma \ref{dual-estimates}
for large values of $x$, we get that for sufficiently large values of   $g^2/\mu(g^2)$
$$
\Phi(u/2) \le F^{-1} \Bigl(1+\frac{u}{2}\Bigr) \le F^{-1} \Bigl(1 + \frac{1}{2} F \Bigl( \frac{g^2}{\mu(g^2)}\Bigr)\Bigr)
\le F^{-1} \Bigl( F \Bigl( \frac{g^2}{\mu(g^2)}\Bigr)\Bigr)
= \frac{g^2}{\mu(g^2)}.
$$
Hence
 $\Phi(u/2)$ is bounded by  $\frac{g^2}{\mu(g^2)} + B$ for
a sufficiently large number $B$
depending only on $F$ and
$\int_{\R^d}(2\Phi(u/2)-1) \,d \mu  \le 2B+1$.
This completes the proof.
\end{proof}

In the following lemma we prove some simple estimates which will be used below.

\begin{lemma}
\label{simple-est}
Suppose that $F$ satisfies assumptions A1)-A3). For
every $K>1$  there exist a number $B$ depending on
$K$ and a number $C$
depending on $K$ and $\Delta$ such that
for every  $f \in L^2(\mu)$ one has
$$
\int_{\{f^2 \ge K \mu(f^2)\}} f^2 F\Bigl(\frac{f^2}{\mu(f^2)}\Bigr) \,d\mu
\le C \cdot \mbox{{\rm Var}}_{\mu} f + \int_{\R^d} f^2 F\Bigl(\frac{f^2}{\mu(f^2)}\Bigr) \,d\mu
$$
$$
\int_{\R^d} f^2 F\Bigl(\frac{f^2}{\mu(f^2)}\Bigr) \,d\mu
\le  B \cdot\mbox{{\rm Var}}_{\mu} f + 2 \int_{\R^d} \bigl( f(x) - \sqrt{K \mu(f^2)}  \bigr)^2_{+}
 F\Bigl(\frac{f^2}{\mu(f^2)}\Bigr) \,d\mu.
$$
\end{lemma}
\begin{proof}
To prove the first estimate we consider
$$
-\int_{\{f^2 \le K \mu(f^2)\}} f^2 F\Bigl(\frac{f^2}{\mu(f^2)}\Bigr) \,d\mu.
$$
Let $\tilde{K} = \min(K,1+ \Delta)$. Since $F(y) \ge 0$ for $y \ge 1$, one has
$$
-\int_{\{f^2 \le K \mu(f^2)\}} f^2 F\Bigl(\frac{f^2}{\mu(f^2)}\Bigr) \,d\mu
\le
-\int_{\{f^2 \le \tilde{K} \mu(f^2)\}} f^2 F\Bigl(\frac{f^2}{\mu(f^2)}\Bigr) \,d\mu.
$$
By the concavity of $-y F(y)$ on $[0,1+\Delta]$ one has
$$
-yF(y) \le (-yF(y))'_{y=1} (y-1) = F'(1) (1-y).
$$
Hence
\begin{align*}
-\int_{\{f^2 \le \tilde{K} \mu(f^2)\}} f^2 F\Bigl(\frac{f^2}{\mu(f^2)}\Bigr) \,d\mu
& \le F'(1)  \mu(f^2) \int_{\{f^2 \le \tilde{K} \mu(f^2)\}} \Bigl( 1 - \frac{f^2}{\mu(f^2)}\Bigr) \,d \mu
\\& =  F'(1)  \mu(f^2) \int_{\{f^2 \ge \tilde{K} \mu(f^2)\}} \Bigl(  \frac{f^2}{\mu(f^2)}-1\Bigr) \,d \mu.
\end{align*}
The desired estimate follows from (\ref{tail-var}).

Let us prove the second estimate. Since $F(y)\ge 0$ for $y \ge K >1 $ and
$$f^2 \le 2K \mu(f^2) + 2 \bigl(f - \sqrt{K \mu(f^2)}\bigr)^2_{+},$$ one has
\begin{align*}
\int_{\R^d} f^2 &  F\Bigl(\frac{f^2}{\mu(f^2)}\Bigr) \,d\mu
\le
\int_{\{f^2 \le K \mu(f^2)\}} f^2  F\Bigl(\frac{f^2}{\mu(f^2)}\Bigr) \,d\mu
+
\int_{\{f^2 \ge K \mu(f^2)\}} f^2  F\Bigl(\frac{f^2}{\mu(f^2)}\Bigr) \,d\mu.
\\&
\le
\int_{\R^d} \min(f^2, K \mu(f^2))  F\Bigl(\frac{f^2}{\mu(f^2)}\Bigr) \,d\mu
+ 2 K \mu(f^2) \int_{\{f^2 \ge K \mu(f^2)\}}   F\Bigl(\frac{f^2}{\mu(f^2)}\Bigr) \,d\mu
\\&
+ 2 \int_{\R^d} \bigl(f - \sqrt{K \mu(f^2)}\bigr)^2_{+}  F\Bigl(\frac{f^2}{\mu(f^2)}\Bigr) \,d\mu.
\end{align*}
The first term on the right-hand side does not exceed
$$
F'(1) \int_{\R^d} \min(f^2, K \mu(f^2))  \Bigl(\frac{f^2}{\mu(f^2)} -1 \Bigr) \,d\mu.
$$
This can be estimated by
$\tilde{C}(K) \cdot \mbox{\rm Var}_{\mu} f$ (see Step 1 in the proof of Theorem \ref{mainth}).
Applying (\ref{tail-var}) and concavity of $F$ we get a similar estimate of the second term. The proof is complete.
\end{proof}

Now we are ready to prove the main result on the tight inequalities.
Following an idea from \cite{GGM} we reduce the problem to $F$-Sobolev  inequalities.
Set $\beta = \frac{\alpha}{\alpha-1}$.
 For every $\tau \ge \frac{2}{\beta}$,
we consider  the following perturbation of $F$:
$$
F_{\tau,\beta } = \psi_{\tau,\beta } (F),
$$
where
$$
\psi_{\tau,\beta}(x) =
\left\{
\begin{array}{lcr}
x, \  \quad        \ x \le 1 \\
 \frac{\beta}{2} \bigl( [1+ \tau(x-1)]^{\frac{2}{\tau \beta}}-1 \bigr) +1, \quad x \ge 1. \\
\end{array}
\right.
$$
 Note that $\psi_{\tau,\beta}(x)$ is a concave increasing function such that $\psi_{\tau,\beta}(x) \le x$.
Obviously,  $\psi_{\frac{2}{\beta},\beta}(x) = x$.

\begin{remark}
\label{pres}
It can be easily verified that this perturbation preserves functions satisfying
 assumptions A1)-A4).
\end{remark}

\begin{theorem}
\label{modtight}
Let $\alpha>1$ and $1 \ge \tau \ge \frac{2}{\beta}$.
Consider the cost function
$c =  c^{*}_{A,\frac{\alpha \tau}{\alpha-1}}$, where $A>0, \varepsilon>0$.
Assume  that $F$, $c$, $\mu$,  and $K$ satisfy the assumptions of Theorem \ref{mainth} for some $K \ge 2$.
 Assume in addition that
\begin{itemize}
\item [1)]
$F$ satisfies assumptions A3)-A4)
\item [2)]
 there exists  $\delta>0$ such that for
 $R = R_{\frac{K-1}{K}}$ one has
$$
 \int_{B^{c}_{R}}\Phi_{\tau, \beta}  \bigl(\delta  |I_{F_{\tau, \beta}}|^2 \bigr) \,d \mu < \infty,
$$
where
$$\Phi_{\tau, \beta}(x)  = \sup_{\{  y > 0\}} \bigl( \bigl<x,y\bigr> - y F_{\tau, \beta}(y) + y \bigr)
= \bigl( yF_{\tau, \beta}(y) - y \bigr)^{*}(x),$$
\item [3)]
$\mu$ satisfies  (\ref{P2}) for some $\lambda_2$.
\end{itemize}
Then  there exist  $B, C >0$ such that the
 following  modified $F$-Sobolev inequality holds:
\begin{align}
\label{Fmod}
&
\int_{\R^d} f^2 F \Bigl( \frac{f^2}{\int_{\R^d} f^2 \,d \mu} \Bigr) \,d \mu
\\&
\nonumber
\le
C\int_{\{f^2 \ge K \int_{\R^d} f^2 \,d \mu\}} f^2 c^{*} \Bigl( \frac{|\nabla f|}{|f|}\Bigr) \,d \mu
+ B \cdot\mbox{\rm Var}_{\mu} f.
\end{align}
\end{theorem}

\begin{proof}
We follow the arguments from \cite{GGM}.
The case $\tau = \frac{2}{\beta}$ follows from  Theorem \ref{Ftight} and Remark \ref{scal}.
Let $\tau > \frac{2}{\beta}$.
Consider a smooth function $f$.  Without loss of generality one can assume that
$\inf_{x \in X} f(x)   > 0$.
If $f$ satisfies the inequality
$$
\int_{\R^d} f^2\,d \mu  \le \frac{1}{2B}
\int_{\R^d} f^2 F \Bigl( \frac{f^2}{\mu(f^2)}
\Bigr) \,d \mu ,
$$
where $B=B(K)$ is the same as
in (\ref{defmodLSI}), then (\ref{Fmod}) follows directly from Theorem \ref{mainth}.
Hence one can assume that
\begin{equation}
\label{ent-l2}
\int_{\R^d} f^2 F \Bigl( \frac{f^2}{\int_{\R^d} f^2 \,d \mu} \Bigr) \,d \mu \le
2B \int_{\R^d} f^2 \,d \mu.
\end{equation}
Note that if $\sup_{x \in X} f^2 \le K \mu(f^2)$, then
by the concavity of $F$
$$
\int_{\R^d} f^2 F \Bigl( \frac{f^2}{\int_{\R^d} f^2 \,d \mu} \Bigr) \,d \mu \le
C(K) {\rm Var}_{\mu} f
$$
 (see the reasoning of Theorem \ref{mainth}, Step 1).
Hence without loss of generality one can assume that there exists $x_0$ such that
$f(x_0) = \sqrt{K \mu(f^2)}$.  Set
$$
g(x) = f(x_0) + (f(x)-f(x_0))_{+} P \Bigl( \frac{f^2}{\mu(f^2)} \Bigr)/P(K),
$$
where
$$
P(x) =  \sqrt{\frac{F(x)}{F_{\tau,\beta}( x)}} =
\sqrt{\frac{F(x)}{\psi_{\tau,\beta}( F(x))}} .
$$
Obviously, $g \ge f$, since $x \mapsto \frac{x}{\psi_{\tau, \beta}(x)}$ is increasing.
In addition, since $\psi_{\tau,\beta}$ is increasing, we get
$$
\psi_{\tau,\beta}\Bigl(F\Bigl(\frac{f^2}{\mu(f^2)}\Bigr)\Bigr) \ge \psi_{\tau,\beta}\bigl(F(K)\bigr)
$$
if $f(x) \ge f(x_0)$. Hence we get by the Cauchy inequality
$$
\int_{\R^d} g^2 \,d \mu \le C_1(K) \Bigl( \int_{\R^d} f^2 \,d \mu
+
\int_{\{f^2 \ge K \mu(f^2) \}} f^2 F \Bigl( \frac{f^2}{\mu(f^2)} \,d \mu \Bigr)
\Bigr)
$$
for some $C_1(K)$.
By Lemma \ref{simple-est}
\begin{equation}
\label{lemma3.3}
\int_{\{f^2 \ge K \mu(f^2)\}} f^2 F \Bigl( \frac{f^2}{\mu(f^2)} \Bigr) \,d \mu
\le \int_{\R^d} f^2 F \Bigl( \frac{f^2}{\mu(f^2)} \Bigr) \,d \mu
+ C_2(K) \mbox{Var}_{\mu} f
.
\end{equation}
Hence by (\ref{ent-l2}) there exists $M=M(K)$ such that
$$\int_{\R^d} g^2 \,d \mu \le M \int_{\R^d} f^2 \,d \mu.
$$
Taking into account that $g \ge f$, one gets
\begin{align*}
\int_{\R^d} \bigl( [g-g(x_0)]_{+}\bigr)^2 F_{\tau,\beta}& \Bigl( \frac{g^2}{\mu(g^2)}\Bigr) \,d \mu
\ge
\int_{\R^d} \bigl( [g-g(x_0)]_{+}\bigr)^2 F_{\tau,\beta} \Bigl( \frac{f^2}{M \mu(f^2)}\Bigr) \,d \mu
\\&
=
\frac{1}{P^2(K)}
\int_{\R^d} \bigl( [f-f(x_0)]_{+}\bigr)^2
F\Bigl( \frac{f^2}{\mu(f^2)}\Bigr)
\frac{F_{\tau,\beta} \Bigl( \frac{f^2}{M \mu(f^2)}\Bigr)}
{F_{\tau,\beta} \Bigl( \frac{f^2}{\mu(f^2)}\Bigr)} \,d \mu.
\end{align*}
By the concavity of $F_{\tau,\beta}$ one has
$\inf_{x \ge 2  M} \frac{F_{\tau,\beta}(x/M)}{F_{\tau,\beta}(x)} = a > 0$.
Hence
\begin{align*}
&
\int_{\R^d} \bigl( [g-g(x_0)]_{+}\bigr)^2 F_{\tau,\beta} \Bigl( \frac{g^2}{\mu(g^2)}\Bigr) \,d \mu
\ge \\&
\frac{a}{P^2(K)}
\int_{\bigl\{ f^2 \ge 2M \mu(f^2) \bigr\}} \bigl( [f-f(x_0)]_{+}\bigr)^2
F\Bigl( \frac{f^2}{\mu(f^2)}\Bigr) \,d \mu
\\& -
\sup_{K \le t \le 2M} \Bigl| \frac{F_{\tau,\beta}(t/M)}{F_{\tau,\beta}(t) }F(t) \Bigr|
\int_{\R^d}  \bigl( [f-f(x_0)]_{+}\bigr)^2 \,d \mu.
\end{align*}
Thus for some $A_1=A_1(K)>0$ one has
\begin{align*}
 \int_{\R^d }  \bigl( [f-f(x_0)]_{+}\bigr)^2
& F\Bigl( \frac{f^2}{\mu(f^2)}\Bigr)\,d \mu
 \le
\\&  A_1 \int_{\R^d} \bigl( [g-g(x_0)]_{+}\bigr)^2
F_{\tau,\beta} \Bigl( \frac{g^2}{\mu(g^2)}\Bigr) \,d \mu
 + A_1
\int_{\R^d}  \bigl( [f-f(x_0)]_{+}\bigr)^2 \,d \mu.
\end{align*}
We observe that the second term on the right-hand side can be estimated by
${\mbox {\rm Var}}_{\mu} f$, since $\bigl(f-f(x_0)\bigr)_{+} \le  \bigl| f - \mu(f) \bigr|$.
 By Lemma \ref{entr-bounds}  we obtain
\begin{align*}
\int_{\R^d} \bigl( [g-g(x_0)]_{+} & \bigr)^2  F_{\tau,\beta}  \Bigl( \frac{g^2}{\mu(g^2)}\Bigr) \,d \mu
\le
\\& 2 \int_{\R^d} \bigl( [g-g(x_0)]_{+}\bigr)^2 F_{\tau,\beta}
\Bigl( \frac{([g-g(x_0)]_{+})^2}{\mu([g-g(x_0)]_{+})^2)}\Bigr)
\,d \mu
+ c' \int_{\R^d}  ([g-g(x_0]_{+})^2 \,d \mu.
\end{align*}

Since
$$\mu(x: g(x) > g(x_0)) = \mu(x: f(x) > f(x_0)) \le \frac{1}{K} \le \frac{1}{2},
$$
$0$ is the median of $(g-g(x_0))_{+}$. Hence by (\ref{P2})
$$
 \int_{\R^d} \bigl( [g-g(x_0)]_{+}\bigr)^2 \,d \mu
 \le \lambda_2 \int_{\R^d} \bigl|\nabla g \bigr|^2 \,d \mu.
$$
By Assumption 2) and  Theorem \ref{Ftight}
$\mu$ satisfies the $F_{\tau,\beta}$-Sobolev inequality, hence
$$
\int_{\R^d} \bigl( [g-g(x_0)]_{+}\bigr)^2 F_{\tau,\beta}
\Bigl( \frac{([g-g(x_0)]_{+})^2}{\mu([g-g(x_0)]_{+})^2}\Bigr) \,d \mu
\le A_2 \int_{\R^d} | \nabla g|^2 \,d \mu.
$$
Combining the estimates obtained above, we get
$$
\int_{\R^d} \bigl( [f-f(x_0)]_{+}\bigr)^2
F\Bigl( \frac{f^2}{\mu(f^2)}\Bigr)\,d \mu
\le C'\Bigl( \int_{\R^d} | \nabla g|^2 \,d \mu
+ {\mbox  {\rm Var} }_{\mu} f
\Bigr).
$$
Let us estimate $\nabla g$. Set
$h = \frac{f^2}{\mu(f^2)}$.
One has
\begin{align*}
\nabla g =
& \\&
 \Bigl[ \frac{\bigl( f-f(x_0)\bigr)_{+} }{P(h) \cdot P(K)} \Bigl( \frac{F'}{\psi_{\tau,\beta}(F)}
- \frac{F\psi_{\tau,\beta}'(F) F'}{\psi_{\tau,\beta}^2(F)} \Bigr)(h)\frac{f}{\mu(f^2)} \Bigr] \nabla f
+
\Bigl[ I_{\{f \ge f(x_0)\}}\frac{P ( h) }{P(K)} \Bigr]  \nabla f.
\end{align*}
Let us show that for some $B_1=B_1(K)>0$ one has
$$|\nabla g|^2 \le B_1 P^2 ( h) |\nabla f|^2.
$$
It is sufficient to verify that
$$
 \frac{\bigl( f-f(x_0)\bigr)_{+} }{P^2(h) } \Bigl( \frac{F'}{\psi_{\tau,\beta}(F)}
- \frac{F\psi_{\tau,\beta}'(F) F'}{\psi_{\tau,\beta}^2(F)} \Bigr)(h)\frac{f }{\mu(f^2)}
$$
is bounded. Since $\frac{f \bigl( f-f(x_0)\bigr)_{+} }{\mu(f^2)}
\le h$ and $P^2 = F/\psi_{\tau,\beta}(F)$, we have to show that
$$
\frac{ x\psi_{\tau,\beta}(F)}{F}  \Bigl( \frac{F'}{\psi_{\tau,\beta}(F)}
- \frac{F\psi_{\tau,\beta}'(F) F'}{\psi_{\tau,\beta}^2(F)} \Bigr)
=
\frac{xF'}{F} - \frac{x\psi_{\tau,\beta}'(F) F'}{\psi_{\tau,\beta}(F)}
=\frac{xF'}{F} \Bigl(1- \frac{F\psi_{\tau,\beta}'(F) }{\psi_{\tau,\beta}(F)}\Bigr)
$$
 is uniformly bounded on $[K,\infty)$. Indeed, it can be verified directly that
$$0 \le \frac{x
\psi_{\tau,\beta}'(x)}{\psi_{\tau,\beta}(x)} \le 1.$$
 The boundedness of $\frac{xF'}{F}$
is obvious. Finally, we obtain
$$
\int_{\R^d} \bigl( [f-f(x_0)]_{+}\bigr)^2
F\bigl(h \bigr)\,d \mu
\le C \int_{\{f^2 \ge K \mu(f^2)\}} | \nabla f|^2   \frac{F(h)}{\psi_{\tau,\beta}
\bigl(F(h) \bigr)} \,d \mu.
$$
The right-hand side can be estimated by
\begin{align*}
C N^{\beta_{\tau}}\int_{\{f^2 \ge K \mu(f^2)\}} f^2  \Bigl|  \frac{\nabla f}{f}  \Bigr|^{\beta_{\tau}} \,d \mu
+
\frac{C}{N^{(\beta_{\tau}/2)^*}} \int_{\{f^2 \ge K \mu(f^2)\}} f^2
 \Bigl| \frac{F(h)}{ \psi_{\tau, \beta}(F(h))}  \Bigr|^{(\beta_{\tau}/2)^*} \,d \mu
\end{align*}
for arbitrary $N$. Here
$$
\beta_{\tau} = \frac{\alpha \tau}{\alpha-1}, \quad \Bigl(\frac{\beta_{\tau}}{2}\Bigr)^{*} =
\frac{\alpha \tau}{2+\alpha (\tau-2)}.
$$
We note that there exists $C' = C'(K)$ such that
for $x \ge K$ one has
$$
\Bigl( \frac{x}{\psi_{\tau,\beta}(x)} \Bigr)^{\frac{\alpha \tau }{2+\alpha(\tau-2) }}
\le C' x^{\bigl(1 - \frac{2(\alpha-1)}{\tau \alpha} \bigr) \bigl( \frac{\alpha \tau}
{2+\alpha (\tau-2)} \bigr) }
= C' x.
$$
Hence  for arbitrary $\varepsilon>0$
and  all sufficiently large $N$ the latter does not
exceed
$$
C N^{\beta_{\tau}}
\int_{\{f^2 \ge K \mu(f^2)\}} f^2  \Bigl|  \frac{\nabla f}{f}  \Bigr|^{\beta_{\tau}} \,d \mu
+ \frac{C C'}{N^{(\beta_{\tau}/2)^*}} \int_{\{f^2 \ge K \mu(f^2)\}}  f^2
F \bigl(h \bigl) \,d \mu.
$$

We recall that
$$
c^{*}(x) =  c_{A,\frac{\alpha \tau}{\alpha-1}}(x) \le \lambda  |x|^{\frac{\alpha \tau}{\alpha-1}}.
$$
for $|x|>1$ and some $\lambda= \lambda(A,\alpha,\tau)$.
Obviously, there exists $a(A,\alpha,K)>0$ such that
$$
|x|^{\beta_{\tau}} \le a(A,\alpha,K) c^{*}(x)
$$
for $x \ge K$.
Hence by (\ref{lemma3.3}) there exists $C=C(\alpha, A, K)$ such that
\begin{align*}
&
\int_{\R^d} \bigl( f-f(x_0)\bigr)^2_{+}
F\Bigl( \frac{f^2}{\mu(f^2)}\Bigr)\,d \mu
\\&
\le C\int_{\{f^2 \ge K \mu(f^2)\}} f^2 c^{*} \Bigl(  \frac{|\nabla f|}{|f|}  \Bigr) \,d \mu
+ C \mbox{Var}_{\mu}  f +
\varepsilon \int_{\R^d}  f^2
F \Bigl( \frac{f^2}{\mu(f^2)}\Bigl) \,d \mu,
\end{align*}
where $\varepsilon$ can be chosen arbitrarily small. It remains
to estimate the last term on the right-hand side by Lemma \ref{simple-est}
$$
\int_{\R^d}  f^2
F \Bigl( \frac{f^2}{\mu(f^2)}\Bigl) \,d \mu \le B(K) \cdot {\mbox {\rm Var}}_{\mu} f
+
2 \int_{\R^d} \bigl( f-\sqrt{K \mu(f^2)} \bigr)^2_{+}
F\Bigl( \frac{f^2}{\mu(f^2)}\Bigr)\,d \mu
$$
and choose a sufficiently small $\varepsilon$.
The proof if complete.
\end{proof}

Now let us apply this result in the case of  a special lower bound for the isoperimetric function.

\begin{theorem}
\label{modtight2}
Let $\varphi$ be a function satisfying assumptions A1)-A4) such that $\varphi(x_0)=1$.
For every $\tau \le 1$ define the corresponding generalized entropy
$$
F(x) = F_{\tau}(x) :=  \left\{
\begin{array}{lcr}
 \varphi(x)  \quad       \mbox{if} \quad 0 < x \le x_0 \\
 \frac{1}{\tau} \bigl( \varphi^{\tau}(x) -1 \bigr) + 1 \quad \mbox{if} \ x \ge x_0. \\
\end{array}
\right.
$$
Assume that
\begin{align}
\label{log-power}
{\mathcal I}_{\mu}(t) \ge C t \varphi \Bigl( \frac{1}{t} \Bigr)^{1-\frac{1}{\alpha}}
\end{align}
for some $1 < \alpha \le 2$ and $ t \le 1/2$. 
Then, whenever
 $1 \ge \tau \ge \frac{2}{\beta} = 2 \Bigl( 1- \frac{1}{\alpha} \Bigr)$
there exists
$C_{\tau} > 0$ such that for every smooth $f$ one has
$$
\int_{\R^d} f^2 F_{\tau} \Bigl( \frac{f^2}{\int_{\R^d} f^2 \,d \mu}\Bigr) \,d \mu
\le
C_{\tau} \int_{\R^d} f^2 c_{A,\frac{\tau \alpha}{\alpha-1}}
\Bigl( \frac{|\nabla f|}{|f|}\Bigr) \,d \mu.
$$
\end{theorem}

\begin{proof}
The result follows from Theorem \ref{modtight}.
Obviously, $F_{\tau}$ satisfies A1)- A4).
Let us show, that $\mu$ satisfies (\ref{P2}). Indeed, it suffices to show that $\mu$ satisfies (\ref{Cheeger}).
But (\ref{Cheeger})  easily follows from (\ref{log-power}), since $\varphi$ is increasing.
Note that $$F_{\tau, \beta} = \psi_{\tau,\beta}(F_{\tau}) =F_{\frac{2}{\beta}} .$$
So it suffices to check that
$$
 \int_{B^{c}_{R_{{(K-1)}/K}}}\Phi_{\tau, \beta}  \bigl(\delta  |I_{F_{\tau, \beta}}|^2 \bigr) \,d \mu < \infty,
 $$
 $$
  \int_{B^{c}_{R_{{(K-1)}/K}}} \Phi_{\tau}  \Bigl(\delta c \bigl(I_{F_{\tau}} \bigr)\Bigr) \,d\mu < \infty.
 $$
 for a sufficiently small number
 $\delta$ and a sufficiently large number
 $K$. Here
 $$\Phi_{\tau} = (y F_{\tau}(y) - y)^{*},$$
 $$
 \Phi_{\tau,\beta} = (y F_{\tau,\beta}(y) - y)^{*}= (y F_{\frac{1}{\beta}}(y) - y)^{*}.
 $$
  Recall that
 the cost function is given by
$$
 c= c^{*}_{A,\bigl( \frac{ \tau \alpha}{\tau-1} \bigr)}
=
 c_{A,\bigl( \frac{\tau \alpha }{\alpha(\tau-1)+1} \bigr)}.
$$
By the
definition of $I_{F_{\tau}}$ and $\mathcal{I}_{\mu}$
for all $r \ge R_{1/2}$ we have
$$
I_{F_{\tau}}(r) \le  \frac{  F_{\tau}(\frac{1}{1-\mu(B_r)})}{ C \varphi^{1-\frac{1}{\alpha} }(\frac{1}{1-\mu(B_r)})}.
$$
Hence
\begin{equation}
\label{IFR}
I_{F_{\tau}}(r) \le c_1 \varphi^{\tau - 1+  \frac{1}{\alpha}}\Bigl(\frac{1}{1-\mu(B_r)}\Bigr)
\end{equation}
and
$$
c(I_{F_{\tau}}(r) ) \le \tilde{c_1}  \varphi^{\tau}\Bigl(\frac{1}{1-\mu(B_r)}\Bigr).
$$
Analogously,
$$
I_{F_{\tau,\beta}}(r) \le  \frac{  F_{\frac{2}{\beta}}(\frac{1}{1-\mu(B_r)})}
{ C \varphi^{1-\frac{1}{\alpha} }(\frac{1}{1-\mu(B_r)})}
\le
c_2 \varphi^{\frac{1}{\beta}}\Bigl(\frac{1}{1-\mu(B_r)}\Bigr)
$$
and
$$
I^2_{F_{\tau,\beta}}(r) \le  \frac{  F_{\frac{2}{\beta}}(\frac{1}{1-\mu(B_r)})}
{ C \varphi^{1-\frac{1}{\alpha} }(\frac{1}{1-\mu(B_r)})}
\le
c^2_2 \varphi^{\frac{2}{\beta}}\Bigl(\frac{1}{1-\mu(B_r)}\Bigr).
$$
Hence by Lemma \ref{dual-estimates}
for some $C_1>0$, $R_0>0$ and sufficiently small $\delta$ one has
$$
\Phi_{\tau}  \Bigl(\delta c \bigl(I_{F_{\tau}}(r) \bigr)\Bigr)
\le
\frac{C_1}{(1-\mu(B_r))^{2 c_1 \delta}}
$$
if $r \ge R_{0}$.
In the same way we obtain
$$
\Phi_{\tau,\beta}  \Bigl(\delta I^2_{F_{\tau,\beta}}(r) )\Bigr)
\le
\frac{C_2}{(1-\mu(B_r))^{2 c_2 \delta}}
$$
if $r \ge R_{0}$.

Hence for a sufficiently small $\delta$ and big $K$ functions
$\Phi_{\tau}  \Bigl(\delta c \bigl(I_{F_{\tau}} \bigr)\Bigr)$ and
$\Phi_{\tau,\beta}  \Bigl(\delta I^2_{F_{\tau,\beta}} )\Bigr)$
are dominated by $\frac{N}{(1-\mu(B_r))^{p}}$ with
 some $p<1$ and $N>0$. Since the mapping
$$x \to 1-\mu(\{y: y \le |x|\})
$$
transforms
$\mu$ into Lebesgue measure on $[0,1]$, we obtain
$$
\int_{B_{R_{(K-1)/{K}}}}\Phi_{\tau}  \Bigl(\delta c \bigl(I_{F_{\tau}}(r) \bigr)\Bigr) \,d \mu
\le N \int_{0}^{1/K} \frac{\,dt}{t^{p}}  < \infty.
$$
The same estimate holds for $\Phi_{\tau,\beta}  \Bigl(\delta I^2_{F_{\tau,\beta}} \Bigr)$.
 Hence assumptions of Theorem \ref{modtight} are fulfilled.
Thus, by Theorem \ref{modtight} we have
$$
\int_{\R^d} f^2 F \Bigl( \frac{f^2}{\int_{\R^d} f^2 \,d \mu} \Bigr) \,d \mu
\le
C\int_{\{f^2 \ge K \int_{\R^d} f^2 \,d \mu\}} f^2 c^{*} \Bigl( \frac{|\nabla f|}{|f|}\Bigr) \,d \mu
+ B \cdot\mbox{\rm Var}_{\mu} f.
$$
By the Poincar{\'e} inequality $\mbox{\rm Var}_{\mu} f \le
\int_{\R^d}| \nabla f|^2 \,d \mu$
 (note that the Poincar{\'e} inequality is fulfilled since (\ref{P2}) holds). One can easily verify that
$|x|^2 \le \tilde{B} c^{*}(x)$ for some $\tilde{B}>0$. Hence
$$
|\nabla f |^2 \le \tilde{B} f^2 c^{*} \Bigl( \frac{|\nabla f|}{|f|} \Bigr)
$$
and
$$
\int_{\R^d} f^2 F \Bigl( \frac{f^2}{\int_{\R^d} f^2 \,d \mu} \Bigr) \,d \mu
\le
( C +\tilde{B} B)
\int_{\R^d} f^2 c^{*} \Bigl( \frac{|\nabla f|}{|f|} \Bigr) \,d\mu.
$$
The proof is complete.
\end{proof}

\section{Application to convex measures}
\begin{lemma}
\label{conv-iso}
Let $\mu$ be a convex measure. Then $\sup_{r \ge R_{1/2}} \frac{I_{\log}(r)}{r} < \infty$.
\end{lemma}
\begin{proof}
We apply the following estimate from \cite{Bobkov}:
\begin{align}
\label{Bob-iso}
2r \mu^{+}(A)
\ge
\mu(A) \log \frac{1}{\mu(A)} + (1-\mu(A)) \log \frac{1}{1-\mu(A)}
+
\log \mu \{|x-x_0| \le r\},
\end{align}
which holds for every convex measure $\mu$, every set $A$, every point $x_0$, and any $r>0$.
Let $\mu(A) \le 1/2 - \varepsilon$, where $\varepsilon >0$. Choose $r$ in such a way that
$\mu(A) = \mu(B^{c}_{r})$.
Then
\begin{equation}
\label{bob-ent}
(1-\mu(A)) \log \frac{1}{1-\mu(A)}
+
\log \mu(B_r)
=
\mu(B^c_r) \log \mu(B_r).
\end{equation}
Pick $\delta = \delta(\varepsilon)$  such  that
$$
\Bigl(\frac{1}{2} + \varepsilon \Bigr)^{\frac{1}{1-\delta}}
\ge
\frac{1}{2} - \varepsilon.
$$
Then
$$\mu(B_r) \ge \Bigl(\frac{1}{2} + \varepsilon\Bigr)
\ge \Bigl( \frac{1}{2} - \varepsilon \Bigr)^{1-\delta} \ge \mu^{1-\delta}(B^c_r).
$$
Therefore,
$$
(1-\delta) \mu(A) \log \frac{1}{\mu(A)} + \mu(B^c_r) \log \mu(B_r) = \mu(A)
\Bigl(\log \frac{\mu(B_r)}{\mu^{1-\delta}(B^c_r)} \Bigr) \ge 0.
$$
Hence by (\ref{bob-ent}) we obtain
$$
(1-\delta) \mu(A) \log \frac{1}{\mu(A)} + (1-\mu(A)) \log \frac{1}{1-\mu(A)}
+
\log \mu \{|x-x_0| \le r\} \ge 0
$$
and
$\frac{2 r}{\delta} \mu^{+}(A) \ge \mu(A) \log \frac{1}{\mu(A)}$. It remains to show that
$$\sup_{R_{1/2} \le r \le R_{1/2+\varepsilon}} \frac{I_{\log}(r)}{r} < \infty.
$$
But this follows easily from
(\ref{Bob-iso}). One has to choose a  sufficiently
large number $\tilde{R}$ such that
$$
\inf_{0 \le \tau \le \varepsilon} (1/2 + \tau) \log \frac{1}{1/2+ \tau}
+ \log \mu(\{x: |x| \le \tilde{R}\}) \ge 0.
$$
Then $I_{\log}(r) \le \tilde{R}$. The proof is complete.
 \end{proof}

\begin{corollary}
\label{isoest}
Let $\mu$ be a convex measure and let
$\varphi$ satisfy assumptions A1) - A4).
Suppose  that $g: \R^{+} \to \R $ is increasing and $\int_{\R^d} e^{g(r)} \,d \mu =1$.
If for some $C>0$ and $1 < \alpha \le 2$ one has
\begin{equation}
\label{12.11.2}
\frac{g(r)}{\varphi^{1-\frac{1}{\alpha}}(e^{g(r)})} \ge C r,
\end{equation}
then
$${\mathcal I}_{\mu}(t) \ge k t \varphi \Bigl( \frac{1}{t}\Bigr)^{1-\frac{1}{\alpha}}
$$
with some $k>0$ and $t \le 1/2$.
\end{corollary}
\begin{proof}
By the previous lemma
$$
\mu^{+}(A)
\ge k_0 \frac{\mu(A) \log \bigl( \frac{1}{\mu(A)}\bigr)}{r}
$$
if $\mu(A) = 1-\mu(B_r)$ and $r \ge R_{1/2}$.
By the Chebyshev inequality
$$
\mu(B^c_r) \le \frac{\int_{\R^d} e^{g(|x|)} \,d\mu(x)}
{e^{g(r)}}= \frac{1}
{e^{g(r)}}.
$$
Hence  by (\ref{12.11.2}) one has
$$
\frac{\log}{\varphi^{1-\frac{1}{\alpha}}} \Bigl( \frac{1}{\mu(B^c_r)}\Bigr)
\ge  Cr
$$
for any $r \ge R_{1/2}$.
Consequently,
$$
\mu^{+}(A)
\ge
 k_0 \frac{\mu(A) \log \bigl( \frac{1}{\mu(A)}\bigr)}{r}
 \ge
 C k_0 \mu(A) \varphi^{1-\frac{1}{\alpha}} \Bigl( \frac{1}{\mu(A)}\Bigr).
$$
The proof is complete.
\end{proof}

{\bf Proof of Theorem \ref{convex}:} Follows from Theorem
\ref{modtight2} and Corollary \ref{isoest}.

\begin{example}
Let $\mu = Z e^{-V} \,dx$ be a convex probability measure on $\R^d$ such that
$V(x) \sim |x| \log^{p}|x|$ with $p>0$ as $|x| \to \infty$.
Suppose that $F$ satisfies A1)-A4) and
$F \sim \log^{\frac{\alpha p}{\alpha-1}} \log |x| $ as $|x| \to \infty$.
Applying Theorem \ref{convex} one gets that for every $A>0$ there exists $C>0$ such that for every
 smooth function $f$ one has
$$
\int_{\R^d} f^2 F\Bigl(\frac{f^2}{\int_{\R^d} f^2 \,d \mu }\Bigr) \,d \mu
\le C \int_{\R^d} f^2 c_{A,\frac{\alpha}{\alpha-1}} \Bigl( \frac{|\nabla f |}{f} \Bigr) \,d \mu.
$$
\end{example}

Finally, we prove an inequality of the type (\ref{Bob-Zeg}).

\begin{theorem}
\label{Bob-Zeg2}
Let $\mu$ be a convex measure such that $\int_{\R^d} e^{\varepsilon |x|^{\alpha}} \,d \mu < \infty$ for
some $1< \alpha $ and $\varepsilon>0$. Then the following inequality holds:
$$
\mbox{\rm Ent}_{\mu} |f|^{\beta} \le C \Bigl[ \int_{\R^d} |\nabla f|^{\beta} \,d\mu
+ \mbox{\rm Var}_{\mu} |f|^{\frac{\beta}{2}}
\bigr].
$$
\end{theorem}
\begin{proof}
Set: $g^2 = |f|^{\beta}$. Apply Theorem \ref{mainth} to $g^2$ in place of $f^2$.
Following the proof we get
$$
\mbox{Ent}_{\mu} |f|^{\beta} \le C \mbox{Var}_{\mu} |f|^{\beta} +
C \int_{|f|^{\beta} \ge K \mu(|f|^{\beta})} I_{\log}\bigl( r_{|f|^{\beta}}(f^{\beta})\bigr)
 |f|^{\beta-1} |\nabla f| \,d \mu
$$
with some $K > 1$.
By the H{\"o}lder inequality for every $\delta>0$ there exists $N(C,\delta)>0$ such that
\begin{align*}
C\int_{|f|^{\beta} \ge K \mu(|f|^{\beta})} I_{\log}\bigl( r_{|f|^{\beta}}(f^{\beta})\bigr)
& |f|^{\beta-1} |\nabla f| \,d \mu
\le \\&
N \int_{\R^d} |\nabla f|^{\beta} \,d\mu +
\delta
\int_{|f|^{\beta} \ge K \mu(|f|^{\beta})} I^{\frac{\beta}{\beta-1}}_{\log}\bigl( r_{|f|^{\beta}}
(f^{\beta})\bigr) |f|^{\beta}
\,d \mu.
\end{align*}
Since $|f| \le C(K, \beta) |f-\mu(f)|$ on $\{|f|^{\beta} \ge K \mu(|f|^{\beta})\}$, we get
by the same arguments as in Theorem \ref{mainth}
\begin{align*}
\delta
\int_{|f|^{\beta} \ge K \mu(|f|^{\beta})} I^{\frac{\beta}{\beta-1}}_{\log}\bigl( r_{|f|^{\beta}}
(f^{\beta})\bigr) |f|^{\beta}
 & \,d \mu
 \\&
\le
\delta C(K, \beta)
\int_{\R^d} I^{\frac{\beta}{\beta-1}}_{\log}\bigl( r_{|f|^{\beta}}(f^{\beta})\bigr) |f - \mu(f)|^{\beta}
  \,d \mu
\\&
\le C
\int_{\R^d}  |f - \mu(f)|^{\beta}  \,d \mu
+ \frac{1}{2}
\mbox{Ent}_{\mu} |f - \mu(f)|^{\beta},
\end{align*}
where $C < \infty$ whenever
$$
\int_{B^{c}_{R}} \exp \Bigl(\delta
I^{\frac{\beta}{\beta-1}}_{\log}\bigl( |x|)\bigr)\Bigr) \,d\mu < \infty
$$
with $R=R_{\frac{K-1}{K}}$.
By Corollary \ref{isoest} one has
$I^{\frac{\beta}{\beta-1}}_{\log}\bigl( |x|\bigr) \le C' |x|^{\frac{\beta}{\beta-1}}
=C' |x|^{\alpha}. $ Hence, choosing a sufficiently
small number $\delta$, we obtain
$$\int_{B^{c}_{R}}
 \exp \Bigl(\delta  I^{\frac{\beta}{\beta-1}}_{\log}\bigl( |x|\bigr)\Bigr) \,d\mu< \infty.$$
  Since
$\mu$ is convex, it satisfies the Cheeger inequality,
hence there exists $C(\beta)$ such that for every $f$
one has
$$
\int_{\R^d}  |f - \mu(f)|^{\beta}  \,d \mu \le C(\beta) \int_{\R^d} |\nabla f|^{\beta} \,d \mu
$$
(see \cite{BH} for the proof).
Finally, we arrive at the estimate
\begin{equation}
\label{2505}
\mbox{Ent}_{\mu} |f|^{\beta} \le C \mbox{Var}_{\mu} |f|^{\frac{\beta}{2}}
+ N' \int_{\R^d} |\nabla f|^{\beta} \,d \mu +
\frac{1}{2} \mbox{Ent}_{\mu} |f-\mu(f)|^{\beta}.
\end{equation}
In particular, applying (\ref{2505}) to $f-\mu(f)$,
we get
\begin{align*}
\mbox{Ent}_{\mu} |f-\mu(f)|^{\beta} \le 2C \int_{\R^d} |f-\mu(f)|^{\beta}
+ & 2N' \int_{\R^d} |\nabla f|^{\beta} \,d \mu
\\& \le \bigl(2CC(\beta) + 2N' \bigr)\int_{\R^d} |\nabla f|^{\beta} \,d \mu.
\end{align*}
Combining this estimate again with (\ref{2505}) we get the claim.
\end{proof}

The author was supported  by the RFBR Grant
04--01--00748 and the  DFG Grant 436 RUS 113/343/0(R).

\end{document}